\documentclass[preprint,11pt]{elsarticle}
\usepackage[dvips]{epsfig}
\usepackage{amscd}
\usepackage{a4wide}
\usepackage{amssymb}
\usepackage{amsthm}
\usepackage{comment}
\usepackage{amsmath}
\usepackage{bm} 
\usepackage{latexsym}
\usepackage[a4paper,inner=2cm,outer=2cm,top=3cm,bottom=3cm]{geometry}
\usepackage{anyfontsize}

\usepackage{blindtext}
\usepackage{upref}
\usepackage{enumerate}
\usepackage{xcolor}
\usepackage{mathtools}

\usepackage[
    colorlinks=true,
    linkcolor=blue,
    citecolor=red,
    urlcolor=blue
]{hyperref}
\hypersetup{
    pdftitle={Inverse problem for a coupled nonlinear Klein-Gordon operator},
    pdfauthor={Parveen Kumar, Soumen Senapati, Manmohan Vashisth}
}

\newtheorem{theorem}{Theorem}[section]
\newtheorem{proposition}[theorem]{Proposition}
\newtheorem{lemma}[theorem]{Lemma}

\newtheorem{remark}[theorem]{Remark}
\newtheorem{corollary}[theorem]{Corollary}   
\numberwithin{figure}{section}
\numberwithin{table}{section}
\allowdisplaybreaks	
\newcommand{\R}{\mathbb{R}}

\renewcommand{\O}{\Omega}

\usepackage{mathrsfs}
\usepackage{tikz}
\usepackage{dsfont}
\usepackage{relsize}
 \mathtoolsset{showonlyrefs}
  \def\diam{ {{\rm diam}}}

\def\supp{ {\rm supp} }

\numberwithin{equation}{section} \allowdisplaybreaks
\usepackage[normalem]{ulem}	
\usepackage{cancel}
\usepackage{cellspace}
\usepackage{makecell}
\begin{document}
\begin{frontmatter}
\title{Inverse Boundary Value Problems for a Coupled Semilinear Klein–Gordon System}
\author[aff1]{Parveen Kumar
}
\ead{parveen.24maz0013@iitrpr.ac.in}

\author[aff2]{Soumen Senapati}
\ead{soumen.senapati@kit.edu}
\author[aff1]{Manmohan Vashisth\corref{cor1}}
\ead{manmohanvashisth@iitrpr.ac.in}

\cortext[cor1]{Corresponding author}

\affiliation[aff1]{
    organization={Department of Mathematics,
    Indian Institute of Technology Ropar},
    city={Rupnagar},
    postcode={140001},
    country={India}
}
\affiliation[aff2]{
    organization={Institut f\"ur Angewandte und Numerische Mathematik, Karlsruher Institut f\"ur Technologie, Englerstr. 2},
    city={Karlsruhe},
    postcode={76131},
    country={Germany}
}
\begin{abstract} 
This article addresses the coefficient identification inverse problem for a coupled nonlinear Klein-Gordon equation with cubic nonlinearity. More precisely, we consider a coupled nonlinear Klein-Gordon equation in a bounded cylinder contained in $\mathbb{R}^{1+\ell}$, where $ \ell\geq 2$, and establish the unique recovery of the time-independent coefficients appearing in both the zeroth-order perturbation and nonlinearity term from both full and partial boundary measurements of the solution. The unique recovery of zeroth-order coefficients is obtained via first-order linearization; however, the nonlinearity coefficients are recovered using the higher-order linearization.
 Finally using the  unique continuation principle, we establish the unique identification of the coefficients from the partial boundary measurements of the solution.
  
\end{abstract}
\end{frontmatter}
\section{Introduction}
\subsection{Formulation of the problem}
For $\ell\geq 2$, let $\Omega\subset\mathbb{R}^\ell$ be an open, bounded, and simply connected set with smooth boundary $\partial \Omega$. Given  $T>0$, we denote by $\O_T:=(0,T)\times \O$ and its lateral boundary by $\partial \O_T:=(0,T)\times\partial\Omega$. In this article, we investigate the inverse problem associated with identifying the coefficients of a coupled nonlinear Klein-Gordon operator given by  
\begin{align}\label{Operator: LPQ}
\mathcal{L}_{\mathbb{P},\mathbb{Q}}\overrightarrow{u}&:=  \begin{bmatrix}
      \partial_t^2 u_1-\Delta u_1 +p_{11}(x) u_1+p_{12}(x)u_2+q_{11}(x)u_1^3+q_{12}(x)u_1^2u_2 \\[2pt]
      \partial_t^2 u_2-\Delta u_2 +p_{21}(x) u_1+p_{22}(x)u_2+q_{21}(x)u_2^2 u_1+q_{22}(x)u_2^3
  \end{bmatrix},
  \end{align}
where 
 \(   \overrightarrow{u}:=(
      u_1, u_2)^\top
\) represents the displacement vector and $\mathbb{P}$, $\mathbb{Q}$ denote the matrix-valued potentials given by 
  \begin{align}
      \mathbb{P}(x):=\begin{bmatrix}
          p_{11}(x) & p_{12}(x)\\
          p_{21}(x)&p_{22}(x)
      \end{bmatrix} \text{ and }   \mathbb{Q}(x):=\begin{bmatrix}
          q_{11}(x) & q_{12}(x)\\
          q_{21}(x)&q_{22}(x)
      \end{bmatrix}.  
  \end{align} 
  Physically, the Klein-Gordon equation describes the dynamics of waves or particles at the subatomic level. The cubic nonlinearity comes into play when quartic scalar fields interact. Moreover, this model is used to study solitary-wave interactions \cite{NSManton1997}, energy transfer \cite{Goodman2005}, and the stability of localized waves \cite{Kopylova2022, Raban2022Stability}.
  
In this article, we consider the following initial boundary value problem (IBVP) associated with the aforementioned Klein-Gordon operator
\begin{equation}\label{Govern Eq}
\begin{cases}
\displaystyle \mathcal{L}_{\mathbb{P},\mathbb{Q}}\overrightarrow{u}(t,x) = \overrightarrow{0},& (t,x)\in \Omega_T,  \\[2pt]
\overrightarrow{u}(t,x)= \overrightarrow{f}(t,x),&  (t,x)\in \partial\Omega_T ,\\[2pt]
\overrightarrow{u}(0,x)= \overrightarrow{0}, \quad\partial_t\overrightarrow{u}(0,x)= \overrightarrow{0}, &\quad x\in \Omega. 
\end{cases}
\end{equation}

We plan to establish the unique recovery of the coefficients $ \mathbb{P},  \mathbb{Q} \in \textbf{C}_c^\infty(\O)$ appearing in \eqref{Govern Eq}, from the knowledge of  Dirichlet-to-Neumann maps measured on whole of  $\partial\O_T$ and a suitable open subset of it. 

Throughout this article,  the boldface notation is used to indicate that the corresponding object is componentwise in the underlying function space. For example, writing $P \in \mathbf{C}_c^\infty(\Omega)$ means that each component of $P$ belongs to $C_c^\infty(\Omega)$. 
In all of our discussion, we use the notation $\Box$ interchangeably to denote the standard wave operator $\partial_t^2-\Delta$. 
In addition, the notation $``\circ"$ is used to denote the component-wise multiplication of two vectors, that is,
  \begin{align}
      \begin{bmatrix}
      v_1\\ v_2
  \end{bmatrix} \circ\begin{bmatrix}
    w_1\\ w_2
  \end{bmatrix}:=\begin{bmatrix}
   v_1 w_1\\v_2 w_2
  \end{bmatrix}.
  \end{align}
We make use of the above notation to rewrite
 $\mathcal{L}_{\mathbb{P},\mathbb{Q}}\overrightarrow{u}$ as 
  \begin{align}
    \mathcal{L}_{\mathbb{P},\mathbb{Q}}\overrightarrow{u}=  (\partial_t^2 -\Delta) ( \mathbb{I}_2 )\overrightarrow{u}+\mathbb{P}\overrightarrow{u}+ (\mathbb{Q}\overrightarrow{u})\circ(\overrightarrow{u}\circ \overrightarrow{u}),
  \end{align}
  where $\mathbb{I}_2$ represents the identity matrix of order $2$. For $\overrightarrow{f} \in \mathtt{F}^\sigma_{m+1}$, where $\mathtt{F}^\sigma_{m+1}$ is given by \eqref{E-delta}, the IBVP \eqref{Govern Eq} admits a unique solution $\overrightarrow{u}\in  \mathbb{B}_R(0)\subset\mathscr{D}_{m+1}$ (we refer the reader to \ref{Well-posedness of the forward model} for the details).
Based on the well-posedness of the forward model \eqref{Govern Eq}, we define the Dirichlet-to-Neumann (DN) map $\Lambda_{\mathbb{P},\mathbb{Q}}(\overrightarrow{f}):\mathtt{F}^\sigma_{m+1} \rightarrow \textbf{H}^{m}(\partial \O_T)$, given  by 
  \begin{align}\label{Def: Full DN map}
      \Lambda_{\mathbb{P},\mathbb{Q}}(\overrightarrow{f}):=\partial_\nu \overrightarrow{u}\rvert_{\partial\O_T},
  \end{align}
  where $\nu$ denotes the unit normal vector pointing outward to $\partial\O$ and $\overrightarrow{u}$ is a solution to IBVP \eqref{Govern Eq}. Set 
  \begin{align}
    \mathbb{P}^k(x):=\begin{bmatrix}
          p_{11}^k(x) & p_{12}^k(x)\\
          p_{21}^k(x)&p_{22}^k(x)
      \end{bmatrix} \text{ and }   \mathbb{Q}^k(x):=\begin{bmatrix}
          q_{11}^k(x) & q_{12}^k(x)\\
          q_{21}^k(x)&q_{22}^k(x)
      \end{bmatrix},  
  \end{align}
  for $k \in \{1,2\}$, where $k$ denotes whether the coefficients belong to $\mathbb{P}^k$ or $\mathbb{Q}^k$; this notation should be interpreted as an index, not an exponentiation operation.
  We now present our first result, which establishes the unique identification of the coefficients from the full DN map and is stated as follows.
   \begin{theorem}
\label{thm:main_result_full_data}
Let $T> \diam(\O)$,  $m>\ell+1$ with $\ell\ge 2$ and  $\mathbb{P}^k,\mathbb{Q}^k \in \textbf{C}^{\infty}_c(\O)$  for $k\in\{1,2\}$. Now if   \[\Lambda_{\mathbb{P}^1,\mathbb{Q}^1}(\overrightarrow{f})=\Lambda_{\mathbb{P}^2,\mathbb{Q}^{2}}(\overrightarrow{f}) \mbox{ for all } \overrightarrow{f} \in \mathtt{F}^\sigma_{m+1},\]then $\mathbb{P}^1=\mathbb{P}^2$, and  $\mathbb{Q}^1=\mathbb{Q}^2$ in  $\Omega$.
\end{theorem}
  Next, we extend the above result and establish the unique identification of the coefficient from the partial DN map $\Lambda^{\sharp}_{\mathbb{P},\mathbb{Q}}(\overrightarrow{f})$, which is defined by  \begin{align}\label{Def: partial DN map}
      \Lambda^{\sharp}_{\mathbb{P},\mathbb{Q}}(\overrightarrow{f}):=\partial_\nu \overrightarrow{u}\rvert_{(0,T)\times\Gamma}, ~~\text{ for all }\overrightarrow{f}\in \mathtt{F}^\sigma_{m+1}, 
  \end{align}
  where
  $\Gamma\subset \partial \Omega$ be a nonempty open set, $T>0$ is sufficiently large and $\overrightarrow{u}$ be the solution to \eqref{Govern Eq}. Now, to formalize the exact statement of the desired problem, we first consider a nonempty  relatively open neighborhood of $\partial\Omega$ in $\overline{\Omega}$, say
  $\mathcal{U}\subseteq \overline{\Omega}$. Also,  we introduce an admissible class $\mathcal{P}$,   given by
\begin{equation}
\label{admissible class}
\begin{aligned}
\mathcal{P}
:= &
\left\{(\mathbb{P},\mathbb{Q})\in \textbf{C}^{\infty}(\overline {\Omega})\times \textbf{C}^{\infty}(\overline {\Omega}): \ 
\mathbb{P}=\mathbb{Q}=0,
\text{ in } \mathcal{U}
\right\}.
\end{aligned}
\end{equation}
Now we describe our next result related to  the unique recovery of coefficients appearing in \eqref{Govern Eq} from the partial DN map and is stated as follows: 
  \begin{theorem}
\label{thm:main_result_partial_data}
Let $T$ is sufficiently large,  $m>\ell+1$ with $\ell\ge 2$ and $\mathbb{P}^k,\mathbb{Q}^k\in \textbf{C}_c^\infty(\O)$. Now if  \[\Lambda^\sharp_{\mathbb{P}^1,\mathbb{Q}^1}(\overrightarrow{f})=\Lambda^\sharp_{\mathbb{P}^2,\mathbb{Q}^{2}}(\overrightarrow{f}) \ \mbox{for every $\overrightarrow{f}\in \mathtt{F}^\sigma_{m+1}$},\] then $\mathbb{P}^1=\mathbb{P}^2$ and $\mathbb{Q}^1=\mathbb{Q}^2$ in $\Omega$, provided   $(\mathbb{P}^1-\mathbb{P}^2,\mathbb{Q}^1-\mathbb{Q}^2)\in \mathcal{P}$.
\end{theorem}
\begin{remark}
A sufficiently large $T$ guarantees the existence of $\diam(\O)<T^{'}<T$, 
used in \S \ref{Recovery of matrix-valued potential from partial DN map} and \S \ref{Recovery of nonlinear potential from partial DN map} to establish the uniqueness of the coefficients given in Theorem \ref{thm:main_result_partial_data}.
    
\end{remark}

\subsection{Existing papers}
 In recent years, the study of inverse problems related to systems of equations has grabbed the attention of researchers. In $1993$, Isakov~\cite{Isakov1993OnUI} introduced the technique of higher-order linearization and solved an inverse problem associated with a nonlinear parabolic partial differential equation (PDE). 
Since then, this methodology has enabled researchers to address inverse problems associated with nonlinear equations. In \cite{Qiu2025UniquenessRF}, the authors considered polynomial nonlinearity in the wave operator and established simultaneous recovery of these coefficients from boundary measurements. 
In \cite{bhardwaj2026reconstructionpotentialdampingcoefficients}, Bhardwaj et al.
considered a coefficient-determination inverse problem for the semilinear damped wave operator with polynomial nonlinearity of order $\ell$, assuming coefficients are time-independent. They derived pointwise reconstruction formulae for the damping coefficient, linear and nonlinear potentials from the DN map.
 For some studies on the scalar nonlinear wave equation, we refer readers to \cite{ Nakamura2021,
NakamuraWatanabeKaltenbacher2009} in the Euclidean setup, and \cite{AcostaUhlmannZhai2022, ChenLassasOksanenPaternain2022, 
feizmohammadi2021inverse,
FeizmohammadiOksanen2022, 
  HintzUhlmannZhai2022b, KaltenbacherLorenzi2007, KaltenbacherRundell2022,
 kian2021determination, Kurylev2018InversePF,
lin2024determining,liu2025partial, oksanen2025inverseproblemconnectionssemilinear, UhlmannZhai2021} in the manifold setup, where the authors used an input–output map, or DN map, or source-to-solution map to determine the coefficients.  For a one-dimensional nonlinear Klein-Gordon equation, we refer to \cite{Tekin2019}, where a time-dependent coefficient multiplied to nonlinear coefficients is uniquely determined using the finite difference scheme.

In comparison to the scalar wave equation, 
the inverse problem associated with the system of equations is studied less. The unique recovery of the time-independent potential in a wave equation has been studied in \cite{avdonin1992boundary}. Eskin and Ralston \cite{eskin2003inverse} uniquely determined both the zeroth-order and first-order matrix-valued coefficients appearing in Schr\"odinger equation from the DN map. Uniqueness of a time-dependent matrix potential is demonstrated in \cite{khanfer2019inverse} in one spatial dimension. 
Waters \cite{Waters2019} established a unique identification of sound speed appearing in a coupled semilinear wave system from the source-to-solution map.
In \cite{mishra2021determining}, Mishra and Vashisth established unique recovery of a zeroth-order matrix-valued potential appearing in a perturbed linear system of wave equations using partial boundary data. Kumar et al. \cite{kumar2024stable} studied an analogous system in an infinite waveguide and derived the stability estimate.
Recently, \cite{BenFraj2025, kumar2026stabilitydeterminationcoefficientscoupled, KumarSarkar2026} discussed the stability in determining the time-dependent matrix potential for the wave equation from partial boundary measurements.

On a global hyperbolic
manifold, Filippas et al. \cite{FilippasOksanenSarkkinen2026} uniquely identified the time-dependent matrix-valued potential in the presence of a connection 1-form (time-independent) using a source-to-solution map. In a Lorentzian manifold, Gomes et al. \cite{Gomes2026} study the unique determination of a connection and potential (up to a natural gauge) appearing in a wave operator
from the Dirichlet-to-Neumann map. For some other works related to the systems, we refer the reader to \cite{
arrepu2026recoverycoefficientsschrodingernonlinear, ChenLassasOksanenEtAl2021, ChenLassasOksanenPaternain2022, CARSTEA2023109740, Hamrouni_2025, Khrabustovskyi2021, Tetlow2022}.

 To the best of the authors knowledge, the present work is the first to deal with inverse problems related to simultaneous recovery of time-independent matrix-valued potentials and nonlinear potentials in a nonlinear coupled wave system from the boundary measurements.

In this article, we uniquely identify the time-independent coefficients from both full and partial DN maps. For the full data setup, we use the geometric optics solutions given in Lemma \ref{lemma: geom opt sol v_{1}} and \ref{lemma: geom opt sol v} to establish the unique recovery of zeroth-order coefficients. However, to establish the unique identification of nonlinear coefficients, we use the asymptotic solutions given in Proposition \ref{Asymptotic solutions for IBVP} and \ref{Asymptotic solutions for backward IBVP}. In the case of the partial DN map, we use the unique continuation argument given in Theorem \ref{UCP} to convert the problem into some sort of full data problem and then use the previous analysis to establish the unique determination of the coefficients.  
The study of the Klein-Gordon equation in \cite{Takei_2023} and the inverse problem discussed in \cite{Bellassoued_Fraj,
bhardwaj2026reconstructionpotentialdampingcoefficients, kumar2026partial} motivates us to consider the problems stated in Theorem \ref{thm:main_result_full_data} and \ref{thm:main_result_partial_data}.
\\

The rest of this article is arranged as follows: \S \ref{Section: Special solutions and unique continuation principle} is devoted to the study of geometric optics solutions. Next, we establish the unique recovery of the zeroth-order potential and nonlinear coefficients from the full DN map in \S \ref{Recovery of the coefficients from full DN map}. In \S \ref{Recovery of the coefficients from partial DN map}, we uniquely determined the coefficients from the partial DN map using the unique continuation principle. Finally, in \ref{Well-posedness of the forward model}, we discussed the well-posedness of the forward model \eqref{Govern Eq}.

\section*{Acknowledgments}
The first author acknowledges financial support from the Council of Scientific and Industrial Research (CSIR), India (File No. 09/1005(19269)/2024-EMR-I). The second author was supported by Deutsche Forschungsgemeinschaft
(DFG, German Research Foundation) -- Project-ID 258734477 -- SFB 1173. 
The last author acknowledges the Anusandhan National Research Foundation (ANRF), Government of India, for financial support provided through the ARG-MATRICS Grant (File No. ANRF/ARGM/2025/002368/MTR).
The authors further acknowledge partial financial support from the Department of Science and Technology (DST), Government of India, under the FIST Programme (Reference No. SR/FST/MS-I/2018/22(C)).

\section{Special solutions}\label{Section: Special solutions and unique continuation principle}
In this section, we recall and extend some existing results needed for the upcoming analysis. We use the term special solutions for geometric optics solutions and asymptotic solutions. We split this section into two subsections: first, we invoke the geometric optics solutions from \cite[Lemma~~(2,3)]{Isakov1991AnIH} and extend them to system of equations in \S \ref{Subsection: Geometric optics solutions}. Next, in \S \ref{Subsection: Asymptotic solutions}, we establish the asymptotic solution associated to initial-value problem \eqref{PDE: linear Forward_problem} and \eqref{PDE: linear Backward_problem}; However, in contrast to geometric optics solutions \eqref{eq: v_{(1)}} and \eqref{eq: v}, the asymptotic solutions are  not directly obtained by multiplying with a constant vector to that one derived for scalar equation \cite{bhardwaj2026reconstructionpotentialdampingcoefficients}, see \S \ref{Subsection: Asymptotic solutions} below.
\subsection{Geometric optics solutions}\label{Subsection: Geometric optics solutions}
 We  begin with introducing the notations used in the forthcoming study. 
Let $\Phi$ solve the transport equation 
\begin{align}\label{transport equation}
    (\partial_t -\theta\cdot \nabla)\Phi(t,x)=0.
\end{align}
In particular, we choose $\Phi$ as follows: For $\varphi\in C_c^{\infty}(\mathbb{R}^\ell)$ and $\theta\in \mathbb{S}^{\ell-1}$, define \begin{align}\label{Def: Phi}
   \Phi:=\varphi(x+t\theta)
\end{align}
 Then $ \Phi$ solves the transport equation \eqref{transport equation}.
Also, let $\mathcal{L}_{\mathbb{P}}$ denotes the linear part of the operator $\mathcal{L}_{\mathbb{P},\mathbb{Q}}$ given in \eqref{Operator: LPQ}, that is,
 \begin{align}\label{Operator: LP}
    \mathcal{L}_{\mathbb{P}}\overrightarrow{u}=  (\partial_t^2 -\Delta) ( \mathbb{I}_2 )\overrightarrow{u}+\mathbb{P}\overrightarrow{u},
  \end{align}
  and $\phi(t,x)=t+x\cdot\theta$. 
Then by  direct computations, we observe that the geometric optics (GO) solution for the  initial value problem (IVP) 
\begin{align}\label{IVP: Lp 1}
     \begin{cases}
         \mathcal{L}_{\mathbb{P}} \overrightarrow{u}(t,x)= \overrightarrow{0}, \ \ &(t,x) \in \Omega_T,\\
     \overrightarrow{u}(0,x)=\partial_t \overrightarrow{u} (0,x) = \overrightarrow{0},  \ \ &\quad x \in\Omega,
     \end{cases}
 \end{align}
 can be derived by extending the GO solution established in \cite[Lemma~~2]{Isakov1991AnIH} for a scalar wave equation, as stated in the following lemma. 
\begin{lemma}\label{lemma: geom opt sol v_{1}}
Let \(\ell \ge 2\) and let \(\Omega \subset \mathbb{R}^\ell\) be a bounded, open, and connected set with smooth boundary \(\partial \Omega\). For any \(\lambda>0\), any \(\varphi \in C_c^\infty(\mathbb{R}^\ell)\) such that \(\operatorname{supp}(\varphi)\cap \Omega=\emptyset\), and any \(\theta \in \mathbb{S}^{\ell-1}\), the IVP  \eqref{IVP: Lp 1}
admits a solution of the form

  \begin{align}\label{eq: v_{(1)}}
  \begin{split}
 \overrightarrow{u}^{(+)}(t,x) = \exp{(\mathrm{i}\lambda \phi)}\Phi(t,x) \overrightarrow{\mathcal{K}}^{(+)}+ \overrightarrow{\mathcal{R}}_{\overrightarrow{\mathcal{K}}^{(+)}}(t,x,\lambda),
 \end{split}
\end{align}
where $\overrightarrow{\mathcal{K}}^{(+)}:=(\mathcal{K}^+_1,\mathcal{K}_2^+)^\top\in \mathbb{R}^2$ is constant vector
and the remainder term $\overrightarrow{\mathcal{R}}_{\overrightarrow{\mathcal{K}}^{(+)}}(t,x,\lambda)$ satisfy the following 
\begin{align}\label{R_{1} initia}
\begin{aligned} 
   \overrightarrow{\mathcal{R}}_{\overrightarrow{\mathcal{K}}^{(+)}}(0,x,\lambda) = \partial_{t}\overrightarrow{\mathcal{R}}_{\overrightarrow{\mathcal{K}}^{(+)}}(0,x,\lambda) =  \overrightarrow{0},   \text{ in } \Omega \ \mbox{and}\\
    \lambda\|\overrightarrow{\mathcal{R}}_{\overrightarrow{\mathcal{K}}^{(+)}}\|_{\textbf{L}^{2}(\Omega_{T})}+\|\partial_{t} \overrightarrow{\mathcal{R}}_{\overrightarrow{\mathcal{K}}^{(+)}}\|_{\textbf{L}^{2}(\Omega_{T})} \leq C \|\varphi\|_{\textbf{H}^{3}{(\R^\ell})}
   \end{aligned} 
\end{align}
for some constant $C>0$, independent of $\lambda$. 
\end{lemma}
Also, the GO solution associated with the backward problem given by 
 \begin{align}\label{IVP: Lp 2}
     \begin{cases}
         \mathcal{L}_{\mathbb{P}^\top} \overrightarrow{u}(t,x)= \overrightarrow{0}, \ \ &(t,x) \in \Omega_T,\\
     \overrightarrow{u}(T,x)=\partial_t \overrightarrow{u} (T,x) = \overrightarrow{0},  \ \ &\quad x \in\Omega,
     \end{cases}
 \end{align}
can  be derived by extending the GO solution established in \cite[Lemma~~3]{Isakov1991AnIH} for a scalar wave equation as stated in the following lemma.
\begin{lemma}\label{lemma: geom opt sol v}
Let \(\ell \ge 2\) and  \(\Omega \subset \mathbb{R}^\ell\) be a bounded, open, and connected set with smooth boundary \(\partial \Omega\). For any \(\lambda>0\), any \(\varphi \in C_c^\infty(\mathbb{R}^\ell)\) such that $\left(\supp(\varphi)\pm T\theta\right)\cap\Omega=\emptyset$, and any \(\theta \in \mathbb{S}^{\ell-1}\), then the backward problem \eqref{IVP: Lp 2} admits a solution of the form 
  \begin{align}\label{eq: v}
  \begin{split}
 \overrightarrow{u}^{(-)}(t,x) = \exp{(-\mathrm{i}\lambda \phi)}\Phi(t,x) \overrightarrow{\mathcal{K}}^{(-)}+ \overrightarrow{\mathcal{R}}_{\overrightarrow{\mathcal{K}}^{(-)}}(t,x,\lambda),
 \end{split}
\end{align}
where $\overrightarrow{\mathcal{K}}^{(-)}:=(\mathcal{K}^-_1,\mathcal{K}^-_2)^\top\in \mathbb{R}^2$ is constant vector
and the remainder term $\overrightarrow{\mathcal{R}}_{\overrightarrow{\mathcal{K}}^{(-)}}(t,x,\lambda)$ satisfying
\begin{align}\label{R_{1} final time}
   \overrightarrow{\mathcal{R}}_{\overrightarrow{\mathcal{K}}^{(-)}}(T,x,\lambda) = \partial_{t}\overrightarrow{\mathcal{R}}_{\overrightarrow{\mathcal{K}}^{(-)}}(T,x,\lambda) =  \overrightarrow{0}   \text{ in } \Omega.
\end{align}
Moreover, $\overrightarrow{\mathcal{R}}_{\overrightarrow{\mathcal{K}}^{(-)}}(t,x,\lambda)$ satisfies the following estimate
\begin{align}\label{eq: estimate of r_{(2)}}
\lambda\|\overrightarrow{\mathcal{R}}_{\overrightarrow{\mathcal{K}}^{(-)}}\|_{\textbf{L}^{2}(\Omega_{T})}+\|\partial_{t} \overrightarrow{\mathcal{R}}_{\overrightarrow{\mathcal{K}}^{(-)}}\|_{\textbf{L}^{2}(\Omega_{T})} \leq C \|\varphi\|_{\textbf{H}^{3}{(\R^\ell})}.
\end{align}
\end{lemma}

    We observe that the unique reconstruction of the matrix-valued potential $\mathbb{P}$ is manageable with the help of GO solutions derived in Lemmas~\ref{lemma: geom opt sol v_{1}} and  \ref{lemma: geom opt sol v}, however, as discussed in \cite[\S3]{bhardwaj2026reconstructionpotentialdampingcoefficients}, the recovery of the nonlinear potential $\mathbb{Q}$  is not feasible due to the products of three or more $\textbf{L}^2$ functions, which need not be dense in $\textbf{L}^1$. This motivates us to look for asymptotic solutions to resolve this issue. 

\subsection{Asymptotic solutions}\label{Subsection: Asymptotic solutions}
In this subsection, we derive the asymptotic  solutions to the IVP given by 
\begin{align}\label{PDE: linear Forward_problem}
     \begin{cases}
         \mathcal{L}_{\mathbb{P}} \overrightarrow{u}(t,x)= \overrightarrow{0}, \ \ &(t,x) \in \Omega_T,\\
     \overrightarrow{u}(0,x)=\partial_t \overrightarrow{u} (0,x) = \overrightarrow{0},  \ \ &\quad x \in\Omega,
     \end{cases}
 \end{align} 
and its formal $L^2-$adjoint given by 
\begin{align}\label{PDE: linear Backward_problem}
     \begin{cases}
         \mathcal{L}_{\mathbb{P}^\top} \overrightarrow{u}(t,x)= \overrightarrow{0}, \ \ &(t,x) \in \Omega_T,\\
     \overrightarrow{u}(T,x)=\partial_t \overrightarrow{u} (T,x) = \overrightarrow{0},  \ \ &\quad x \in\Omega. 
     \end{cases}
 \end{align}
 These solutions are generalization of  the asymptotic  solutions established  in \cite[\S 3]{bhardwaj2026reconstructionpotentialdampingcoefficients} for the scalar equations. In contrast to the geometric optics result above, simply multiplying the constant vector by the one derived in \cite{bhardwaj2026reconstructionpotentialdampingcoefficients} for the scalar equation doesn't work here. 
So, we begin by following  ansatz  
    \begin{equation}\label{ansatz}
 \overrightarrow{u}(t,x)= \exp{(\mathrm{i}\lambda \phi)}\left(z_{0}(t,x)\overrightarrow{K}+\frac{Z_{1,\overrightarrow{K}}(t,x)}{\lambda}+\frac{Z_{2,\overrightarrow{K}}(t,x)}{\lambda^2}+\cdots+\frac{Z_{N,\overrightarrow{K}}(t,x)}{\lambda^N}\right)+\overrightarrow{\mathcal{R}}_{\overrightarrow{K}}(t,x,\lambda),
    \end{equation}
    as solution to the IVP \eqref{PDE: linear Forward_problem}, 
    where $\overrightarrow{K}\in \mathbb{R}^2$ is constant vector, $Z_{j,\overrightarrow{K}}, \overrightarrow{\mathcal{R}}_{\overrightarrow{K}} $ are $2\times 1$ column vectors for $j=1, \cdots, N$ and $\lambda>1$ is a real parameter.
    Let us denote \[\overrightarrow{U}=\begin{bmatrix}
        U_1\\U_2
\end{bmatrix}:=z_{0}\overrightarrow{K}+\frac{Z_{1,\overrightarrow{K}}}{\lambda}+\frac{Z_{2,\overrightarrow{K}}}{\lambda^2}+\cdots+\frac{Z_{N,\overrightarrow{K}}}{\lambda^N}.\]
    Through a direct computation, one can see that operating $\mathcal{L}_\mathbb{P}$ on above ansatz yields
\begin{align}\label{operator: Lp after substitution ansatz}
\overrightarrow{0}=\mathcal{L}_\mathbb{P}\overrightarrow{u}&=\mathcal{L}_\mathbb{P}(\exp{(\mathrm{i}\lambda \phi)}\overrightarrow{U}+\overrightarrow{\mathcal{R}}_{\overrightarrow{K}})\\&= \exp{(\mathrm{i}\lambda \phi)}\begin{bmatrix}
   \Box U_1+2\mathrm{i}\lambda \partial_tU_1-2\mathrm{i}\lambda \theta\cdot \nabla U_1+p_{11}U_1+p_{12}U_2\\ \Box U_2+2\mathrm{i}\lambda \partial_tU_2-2\mathrm{i}\lambda \theta\cdot \nabla U_2+p_{21}U_1+p_{22}U_2
\end{bmatrix}+\mathcal{L}_\mathbb{P} \overrightarrow{\mathcal{R}}_{\overrightarrow{K}}.
    \end{align}

We now choose $z_0 \overrightarrow{K}:=Z_0$ and $Z_{j,\overrightarrow{K}}$  be the solution to the transport equations 
\begin{align}\label{z_0 transport}
    \begin{cases}
       \mathcal{T}z_0(t,x):= \left(\partial_t - \theta\cdot\nabla\right)z_{0}(t,x) = 0, \ \ &(t,x)\in  \Omega_T,\\
         z_0(0,x)=\partial_t z_0(0,x)=0 \ \ &\quad  x\in  \Omega,
    \end{cases}
\end{align}
 and
\begin{align}\label{remaining transport}
    \begin{aligned}
        \begin{cases}
             \mathcal{T} Z_{j,\overrightarrow{K}}(t,x)=- \dfrac{1}{ 2\mathrm{i}} \mathcal{L}_{\mathbb{P}} Z_{j-1,\overrightarrow{K}}(t,x),\ \ &(t,x)\in  \Omega_T,\\Z_{j,\overrightarrow{K}}(0,x)=\partial_t Z_{j,\overrightarrow{K}}(0,x)=\overrightarrow{0},\ \ &\quad x\in  \Omega,
        \end{cases}
    \end{aligned}
\end{align}
for each $j = 1,\cdots, N$. Note that the solution of the above transport equations exist and are given by 
\begin{align}\label{Def: z0}
   z_0(t,x):=\varphi(x+t\theta),
\end{align}
 and 
 \begin{align}\label{eqn:mj-explicit}
Z_{j,\overrightarrow{K}}(t,x)
=
-\frac{1}{2\mathrm{i}}\int_{0}^{t}
 \mathcal{L}_{\mathbb{P}} Z_{j-1,\overrightarrow{K}}(s,x+(t-s)\theta) 
\,ds ,
\end{align} respectively, where $\varphi\in C_c^{\infty}(\mathbb{R}^\ell)$ with $\supp \varphi \cap \O=\emptyset$, and $\theta\in \mathbb{S}^{\ell-1}$.
  From the above transport equations, we observe that all $\lambda^{-k}$ terms vanish and $\overrightarrow{\mathcal{R}}_{\overrightarrow{K}}(t,x,\lambda)$ solves
\begin{align}
  \mathcal{L}_\mathbb{P} \overrightarrow{\mathcal{R}}_{\overrightarrow{K}}(t,x) =-\lambda^{-N}\exp{(\mathrm{i}\lambda \phi)}\mathcal{L}_\mathbb{P} Z_{N,\overrightarrow{K}}(t,x).
\end{align}
Since we require that $\overrightarrow{u}$ defined by \eqref{ansatz} solves the IVP \eqref{PDE: linear Forward_problem}, the problem reduces to establishing the existence of $\overrightarrow{\mathcal{R}}_{\overrightarrow{K}}$ satisfying
\begin{align}
\begin{cases}
  \mathcal{L}_\mathbb{P} \overrightarrow{\mathcal{R}}_{\overrightarrow{K}}(t,x) =-\lambda^{-N}\exp{(\mathrm{i}\lambda \phi)}\mathcal{L}_\mathbb{P} Z_{N,\overrightarrow{K}}(t,x),\ \ &(t,x)\in  \Omega_T, \\
   \overrightarrow{\mathcal{R}}_{\overrightarrow{K}}(0,x)= \partial_t \overrightarrow{\mathcal{R}}_{\overrightarrow{K}}(0,x)=\overrightarrow{0}, \ &\quad  x\in  \Omega.
\end{cases}
\end{align} 
Using the well-posedness result given by Proposition \ref{proposition 1 system} and Sobolev embedding, we have $\overrightarrow{\mathcal{R}}_{\overrightarrow{K}}\in \textbf{C}^\infty(\O_T)$ satisfying 
\begin{align}
    \lVert \overrightarrow{\mathcal{R}}_{\overrightarrow{K}} \rVert_{\textbf{H}^m(\O_T)}\leq C\lambda^{m-N}.
\end{align}
Again, we use the Sobolev embedding theorem $(m>\ell+1)$ to conclude that
\begin{align}
   \lVert \overrightarrow{\mathcal{R}}_{\overrightarrow{K}} \rVert_{\textbf{L}^\infty(\O_T)}\leq  C \lVert \overrightarrow{\mathcal{R}}_{\overrightarrow{K}} \rVert_{\textbf{H}^m(\O_T)}\leq C\lambda^{m-N}.
\end{align}
For $N\geq m+1$, we have
\begin{align}
   \lVert \overrightarrow{\mathcal{R}}_{\overrightarrow{K}} \rVert_{\textbf{L}^\infty(\O_T)}\leq C\lambda^{-1},
\end{align}
for any $\lambda>1$.
 On the basis of the preceding analysis, we have established  the following proposition.
\begin{proposition}[Asymptotic solutions for forward problem]\label{Asymptotic solutions for IBVP}
 Let $\mathbb{P} \in \textbf{C}_c^{\infty}(\Omega)$ and $T>0$ be fixed. Then, the IVP \eqref{PDE: linear Forward_problem}
has an asymptotic solution  of the form 
    \begin{equation}\label{eqn: IVP geom opt sol 1}
 \overrightarrow{u}(t,x)= \exp{(\mathrm{i}\lambda \phi)}\left(z_{0}(t,x)\overrightarrow{K}+\frac{Z_{1,\overrightarrow{K}}(t,x)}{\lambda}+\frac{Z_{2,\overrightarrow{K}}(t,x)}{\lambda^2}+\cdots+\frac{Z_{N,\overrightarrow{K}}(t,x)}{\lambda^N}\right)+\overrightarrow{\mathcal{R}}_{\overrightarrow{K}}(t,x,\lambda),
    \end{equation}
    where $\overrightarrow{K}\in \mathbb{R}^2$ is constant vector, $z_0$ solves the  transport equation
\begin{align}\label{IVP: transport equation z_0}
\begin{cases}
     \mathcal{T}z_0(t,x):= \left(\partial_t - \theta\cdot\nabla\right)z_{0}(t,x) = 0, \ \ &(t,x)\in  \Omega_T, \\[2pt]
    z_0(0,x)=\partial_t z_0(0,x)=0, \ \  &\quad x\in  \Omega ,
\end{cases}
\end{align}
  and the finite sequence of vectors $ Z_{j,\overrightarrow{K}}$ solves the system of transport equations:
\begin{align}\label{eq:transport--mj}
\begin{cases}
  \mathcal{T} Z_{j,\overrightarrow{K}}(t,x)=-\frac{1}{  2\mathrm{i}}  \mathcal{L}_{\mathbb{P}} Z_{j-1,\overrightarrow{K}}(t,x),\ \ &(t,x)\in  \Omega_T, \\[2pt]
    Z_{j,\overrightarrow{K}}(0,x)= \partial_t Z_{j,\overrightarrow{K}}(0,x)=\overrightarrow{0}, \ &\quad  x\in  \Omega,
\end{cases}
\end{align}
for $j = 1,\dots, N$ and $Z_0:=z_0\overrightarrow{K}$. Also, the
correction vector $\overrightarrow{\mathcal{R}}_{\overrightarrow{K}}(t,x,\lambda)$ solves
\begin{align}\label{eqn: reminder}
\begin{cases}
  \mathcal{L}_\mathbb{P} \overrightarrow{\mathcal{R}}_{\overrightarrow{K}}(t,x) =-\lambda^{-N}\exp{(\mathrm{i}\lambda \phi)}\mathcal{L}_\mathbb{P} Z_{N,\overrightarrow{K}}(t,x),\ \ &(t,x)\in  \Omega_T, \\[2pt]
\overrightarrow{\mathcal{R}}_{\overrightarrow{K}}(0,x)= \partial_t \overrightarrow{\mathcal{R}}_{\overrightarrow{K}}(0,x)=\overrightarrow{0}, \ &\quad  x\in  \Omega,
\end{cases}
\end{align}
and satisfies
\begin{align}\label{eq: r__{(2)}}
\|\overrightarrow{\mathcal{R}}_{\overrightarrow{K}}\|_{\textbf{L}^{\infty}(\Omega_{T})}\leq C\lambda^{-1}, \ \text{for any $\lambda>1$},
\end{align}
where $C>0$ is independent of $\lambda$ and $N\geq m+1$.
    \end{proposition}
Using a similar approach, one can derive the asymptotic solutions for the backward problem, as stated in the following proposition.
\begin{proposition}[Asymptotic solutions for backward problem]\label{Asymptotic solutions for backward IBVP}
 Let $\mathbb{P} \in \textbf{C}_c^{\infty}(\Omega)$ and $T>0$ be fixed. Then, the backward problem \eqref{PDE: linear Backward_problem}
has an asymptotic solution  of the form 
    \begin{align}\label{eqn: ansatz2}
 \overrightarrow{u}(t,x)= \exp{(-3\mathrm{i}\lambda \phi)}\bigg(w_{0}\overrightarrow{K}^{(-)}+\frac{W_{1,\overrightarrow{K}^{(-)}}}{\lambda}+\frac{W_{2,\overrightarrow{K}^{(-)}}}{\lambda^2}+\cdots+\frac{W_{N,\overrightarrow{K}^{(-)}}}{\lambda^N}\bigg)+\overrightarrow{\mathcal{R}}_{\overrightarrow{K}^{(-)}}(t,x,\lambda),
    \end{align}
    where $\overrightarrow{K}^{(-)}\in \mathbb{R}^2$ is constant vector and  $W_{j,\overrightarrow{K}^{(-)}},\overrightarrow{\mathcal{R}}_{\overrightarrow{K}^{(-)}}$ are $2\times 1$ column vectors for $j=1,\cdots, N$. Here, $w_0$ solves the  transport equation
\begin{align}\label{IVP: transport w_0}
\begin{cases}
     \mathcal{T}w_0(t,x):= \left(\partial_t - \theta\cdot\nabla\right)w_{0}(t,x) = 0, \ \ &(t,x)\in  \Omega_T, \\
    w_0(T,x)=\partial_t w_0(T,x)=0, \ \  &\quad x\in  \Omega ,
\end{cases}
\end{align}
  and the finite sequence of vectors $ W_{j,\overrightarrow{K}^{(-)}}$ solves the system of transport equations:
\begin{align}\label{eq:transport--Wj}
\begin{cases}
 \mathcal{T} W_{j,\overrightarrow{K}^{(-)}}(t,x)=\frac{1}{6\mathrm{i}}  \mathcal{L}_{\mathbb{P}^\top} W_{j-1,\overrightarrow{K}^{(-)}}(t,x),\ \ &(t,x)\in  \Omega_T, \\
    W_{j,\overrightarrow{K}^{(-)}}(T,x)= \partial_t W_{j,\overrightarrow{K}^{(-)}}(T,x)=\overrightarrow{0}, \ &\quad  x\in  \Omega,
\end{cases}
\end{align}
for $j = 1,\dots, N$ and $W_0:=w_0\overrightarrow{K}^{(-)}$. Also, the
correction vector $\overrightarrow{\mathcal{R}}_{\overrightarrow{K}^{(-)}}(t,x,\lambda)$ solves
\begin{align}\label{eqn: reminder2}
\begin{cases}
  \mathcal{L}_{\mathbb{P}^\top} \overrightarrow{\mathcal{R}}_{\overrightarrow{K}^{(-)}}(t,x) =-\lambda^{-N}\exp{(-3\mathrm{i}\lambda \phi)}\mathcal{L}_{\mathbb{P}^\top} W_{N,\overrightarrow{K}^{(-)}}(t,x),\ \ &(t,x)\in  \Omega_T, \\
   \overrightarrow{\mathcal{R}}_{\overrightarrow{K}^{(-)}}(T,x)= \partial_t \overrightarrow{\mathcal{R}}_{\overrightarrow{K}^{(-)}}(T,x)=\overrightarrow{0}, \ &\quad  x\in  \Omega,
\end{cases}
\end{align}
and satisfies
\begin{align}\label{eq: reminder2 bound}
\|\overrightarrow{\mathcal{R}}_{\overrightarrow{K}^{(-)}}\|_{\textbf{L}^{\infty}(\Omega_{T})}\leq C\lambda^{-1}, \ \text{for any $\lambda>1$},
\end{align}
where $C>0$ is independent of $\lambda$ and $N\geq m+1$.
    \end{proposition}

\section{Proof of Theorem \ref{thm:main_result_full_data}}\label{Recovery of the coefficients from full DN map}
This section is devoted to proving the unique recovery of the coefficients (both the zeroth-order potential and the nonlinear potential) from the full DN map. We divide this section into two subsections, where \S \ref{Recovery of zeroth-order coefficients} and \S\ref{Recovery of the coefficients of nonlinearity} address the unique identification of the zeroth-order and nonlinear coefficients, respectively.
\subsection{Unique determination of matrix-valued potential}\label{Recovery of zeroth-order coefficients} In the following, we start by linearizing the problem using a first-order linearization method. This technique is commonly practiced in the study of nonlinear equations; see, for instance, \cite{bhardwaj2026reconstructionpotentialdampingcoefficients, CARSTEA2019121, kumar2026partial, Lassas2021,    NakamuraWatanabeKaltenbacher2009}.
Given  $f_{11},f_{12},f_{13}, f_{21}, f_{22},f_{23} \in H^{m+1}(\partial \O_T)$, we choose $\varepsilon:=(\varepsilon_1,\varepsilon_2,\varepsilon_3)$, where $\varepsilon_i\geq0$ for all $i=1,2,3$, so that 
\[\varepsilon\overrightarrow{f}:=\begin{bmatrix}
    \varepsilon_1f_{11}+\varepsilon_2f_{12}+\varepsilon_3f_{13}
\\\varepsilon_1f_{21}+\varepsilon_2f_{22}+\varepsilon_3f_{23}\end{bmatrix} \in \mathtt{F}^\sigma_{m+1},\]
for $\sigma>0$ sufficiently small.
Let $\overrightarrow{u}_{\varepsilon\overrightarrow{f}}:=\overrightarrow{u}$ be the solution to the IBVP
\begin{equation}\label{govern eqn: For IBVP in epsilon form}
\left\{\begin{array}{r l c}
\displaystyle \mathcal{L}_{\mathbb{P},\mathbb{Q}}\overrightarrow{u}(t,x) &= \overrightarrow{0},\quad& (t,x)\in \Omega_T,  \\[2pt]
\overrightarrow{u}(t,x)&= \varepsilon\overrightarrow{f}(t,x),\quad&  (t,x)\in \partial\Omega_T ,\\[2pt]
\overrightarrow{u}(0,x)&= \overrightarrow{0}, \quad\partial_t\overrightarrow{u}(0,x)= \overrightarrow{0}, \quad& x\in \Omega.
\end{array}\right.
\end{equation}
The corresponding DN map is given by 
  \begin{align}\label{Def: DN map for IBVP of epsilon form}
      \Lambda_{\mathbb{P},\mathbb{Q}}(\varepsilon\overrightarrow{f}):=\partial_\nu \overrightarrow{u}\rvert_{\partial\O_T}.
  \end{align}
For $i\in\{1,2,3\}$, we apply $\partial_{\varepsilon_i}$ to \eqref{govern eqn: For IBVP in epsilon form} and evaluate at $\varepsilon=0$ to get
\begin{equation}\label{PDE: v power i}
\left\{\begin{array}{r l  c}
\displaystyle
\mathcal{L}_{\mathbb{P}}
\overrightarrow{v}^{(i)}(t,x)&=\overrightarrow{0},\quad& (t,x)\in \Omega_T,  \\[2pt]
\overrightarrow{v}^{(i)}(t,x)&=(f_{1i}, f_{2i})^\top(t,x),\quad&  (t,x)\in \partial\Omega_T ,\\[2pt]
\overrightarrow{v}^{(i)}(0,x)&= \partial_t \overrightarrow{v}^{(i)}(0,x)= \overrightarrow{0}, \quad& x\in \Omega,
\end{array}\right.
\end{equation}
where $\partial_{\varepsilon_i} \overrightarrow{u}\rvert_{\varepsilon=0}:=\overrightarrow{v}^{(i)}=(v_1^{(i)}, v_2^{(i)})^\top$. Also, the first-order linearization of the DN map \eqref{Def: DN map for IBVP of epsilon form} becomes
\begin{align}\label{Def: DN map first-order linearization}
{}^1\Lambda_\mathbb{P}(f_{1i},f_{2i})\lvert_{\partial\O_T}:=\left(\partial_{\varepsilon_i}\Lambda_{\mathbb{P},\mathbb{Q}}(\varepsilon\overrightarrow{f})\lvert_{\varepsilon=0}\right)\Big\lvert_{\partial\O_T}=\partial_\nu \overrightarrow{v}^{(i)}\lvert_{\partial\O_T}.
\end{align}
To avoid the cumbersomeness, we denote $\overrightarrow{v}^{(i)}:=\overrightarrow{v}$, $f_{1i}:=f$, and  $f_{2i}:=g$ throughout the establishment of unique recovery of zeroth-order potential. As a result, \eqref{PDE: v power i} and \eqref{Def: DN map first-order linearization} transformed into
\begin{equation}\label{PDE: first order linearization in simplified form}
\left\{\begin{array}{r l  c}
\displaystyle \mathcal{L}_{\mathbb{P}}\overrightarrow{v}(t,x)&=\overrightarrow{0},\quad& (t,x)\in \Omega_T, \\[2pt]
\overrightarrow{v}(t,x)&= (f,g)^\top(t,x),\quad&  (t,x)\in \partial\Omega_T ,\\[2pt]
\overrightarrow{v}(0,x)&= \partial_t \overrightarrow{v}(0,x)=\overrightarrow{0}, \quad& x\in \Omega,
\end{array}\right.
\end{equation}
 and \begin{align}\label{Def: DN map first order in simplified form}
{}^1\Lambda_\mathbb{P}(f,g)\lvert_{\partial\O_T}=\partial_\nu \overrightarrow{v}\lvert_{\partial\O_T}, \quad f,g\in H^{m+1}(\partial\O_T),
 \end{align}
 respectively.
Next, we assert that the matrix $\mathbb{P}$ can be uniquely determined from the boundary measurement \eqref{Def: DN map first order in simplified form}. To begin,  we denote $\overrightarrow{v}^{j}$ be the solution to the IBVP   
\begin{equation}\label{PDE: compact form vi}
\left\{\begin{array}{r l  c}
\displaystyle \mathcal{L}_{\mathbb{P}^j}
\overrightarrow{v}^j(t,x)& = \overrightarrow{0},\quad& (t,x)\in \Omega_T,  \\[2pt]
\overrightarrow{v}^j(t,x)&= (f, g)^\top(t,x),\quad&  (t,x)\in \partial\Omega_T ,\\[2pt]
\overrightarrow{v}^j(0,x)&= \partial_t \overrightarrow{v}^j(0,x)=\overrightarrow{0}, \quad& x\in \Omega.
\end{array}\right.
\end{equation}
We denote $\overrightarrow{v}^{12}:=\overrightarrow{v}^1-\overrightarrow{v}^2$, and $\mathbb{P}^{12}:=\mathbb{P}^1-\mathbb{P}^2$, where for each $j=1,2$, $\overrightarrow{v}^j$ is solution to \eqref{PDE: compact form vi} and
\begin{align}
   \mathbb{P}^j:= \begin{bmatrix}
        p_{11}^j&p_{12}^j\\[2pt]p_{21}^j
&p_{22}^j    \end{bmatrix}.
\end{align}
Thus, we have 
\begin{align}\label{PDE: used in integral identity}
    \left\{\begin{array}{l  c}
\displaystyle 
\Box \overrightarrow{v}^{12} +\mathbb{P}^1\overrightarrow{v}^{12} = -\mathbb{P}^{12}\overrightarrow{v}^2,\quad& (t,x)\in \Omega_T,  \\[2pt]
\overrightarrow{v}^{12}(t,x)= \overrightarrow{0},\quad&  (t,x)\in \partial\Omega_T ,\\[2pt]
\overrightarrow{v}^{12}(0,x)= \partial_t \overrightarrow{v}^{12}(0,x)=\overrightarrow{0}, \quad& x\in \Omega.
\end{array}\right.
\end{align}
Since \(
      \Lambda_{\mathbb{P}^1,\mathbb{Q}^1}(\varepsilon\overrightarrow{f})= \Lambda_{\mathbb{P}^2,\mathbb{Q}^2}(\varepsilon\overrightarrow{f}),
\) the first order linearization then leads us to
\[{}^1\Lambda_{\mathbb{P}^1}(f,g)\lvert_{\partial\O_T}={}^1\Lambda_{\mathbb{P}^2}(f,g)\lvert_{\partial\O_T}, \quad f,g\in H^{m+1}(\partial\O_T).\]
Consequently, we have $\partial_\nu \overrightarrow{v}^{12}\lvert_{\partial\O_T}=\overrightarrow{0}$.
Next, we
consider $\overrightarrow{v_0}:=(v_{01},v_{02})^\top$ to be the solution to the backward problem 
\begin{align}\label{PDE: backward problem used in integral identity}
    \begin{aligned}
      \left\{\begin{array}{l  c}
\displaystyle 
\Box \overrightarrow{v_0} +(\mathbb{P}^1)^\top\overrightarrow{v_0} = \overrightarrow{0},\quad& (t,x)\in \Omega_T,  \\[2pt]
\overrightarrow{v_0}(T,x)= \partial_t \overrightarrow{v_0}(T,x)=\overrightarrow{0}, \quad& x\in \Omega.
\end{array}\right.
    \end{aligned}
\end{align}
To derive the desired integral identity, we then  multiply \eqref{PDE: first order linearization in simplified form} by $\overrightarrow{v_0}$ and integrate over $\O_T$, which leads us to
\begin{align}
    \int_{\O_T} \left(\Box \overrightarrow{v}^{12} +\mathbb{P}^1\overrightarrow{v}^{12}\right)\cdot \overrightarrow{v_0}~dxdt = -    \int_{\O_T}\mathbb{P}^{12}\overrightarrow{v}^2 \cdot \overrightarrow{v_0}~dxdt
\end{align}
Armed with the initial and final time information available in \eqref{PDE: used in integral identity} and \eqref{PDE: backward problem used in integral identity} respectively, 
the integration by parts  formula yields
\begin{align}\label{integral identity 1}
    \begin{aligned}
           \int_{\O_T}\mathbb{P}^{12}\overrightarrow{v}^2 \cdot \overrightarrow{v_0}~dxdt=0.
    \end{aligned}
\end{align}
We are now in position to utilize the GO solutions established in \S \ref{Subsection: Geometric optics solutions} to conclude about the unique identification of the coefficients. Although to work in a common class of functions $\varphi$ satisfying the transport equation \eqref{transport equation} with $\supp(\varphi)\cap\Omega=\emptyset$ and $ (\supp(\varphi)\pm T\theta)\cap \Omega=\emptyset$, we introduce a function space $\mathscr{C}_c^{\infty}(\mathbb{R}^\ell)$, a subspace of  $ C_c^{\infty}(\mathbb{R}^\ell)$, and it is given by
\begin{align}\label{eq; definition of C(Rn)}
\mathscr{C}_c^{\infty}(\mathbb{R}^\ell):=\{\varphi\in C_c^{\infty}(\mathbb{R}^\ell):\supp(\varphi)\cap\Omega=\emptyset\ \mbox{and}\ (\supp(\varphi)\pm T\theta)\cap \Omega=\emptyset,\ \forall \ \theta\in\mathbb{S}^{\ell-1}\}.
\end{align}
For $\mathcal{K}^{\pm}:=(\mathcal{K}^{\pm}_{1},\mathcal{K}^{\pm}_{2})^{\top}\in\mathbb{R}^2$,  $\varphi\in \mathscr{C}_c^{\infty}(\mathbb{R}^\ell)$, choose (see Lemma \ref{lemma: geom opt sol v_{1}} and Lemma \ref{lemma: geom opt sol v})
\begin{align}\label{GO solution v2}
    \overrightarrow{v}^{2}(t,x) = \varphi(x+t\theta)\exp{(\mathrm{i}\lambda \phi)}\overrightarrow{\mathcal{K}}^{(+)} + \overrightarrow{\mathcal{R}}^{v_2}_{\overrightarrow{\mathcal{K}}^{(+)}}(t,x,\lambda),
\end{align}
and 
\begin{align}\label{GO solution y}
     \overrightarrow{v_0}(t,x) = \varphi(x+t\theta)\exp{(-\mathrm{i}\lambda \phi)}\overrightarrow{\mathcal{K}}^{(-)} + \overrightarrow{\mathcal{R}}^{v_0}_{\overrightarrow{\mathcal{K}}^{(-)}}(t,x,\lambda),
\end{align}
the solutions to IVP \eqref{IVP: Lp 1} and backward problem \eqref{IVP: Lp 2}. Now we substitute \eqref{GO solution v2} and \eqref{GO solution y} in  \eqref{integral identity 1} to arrive at
\begin{align}\label{integral identity opened form}
    \begin{aligned}
        \int_{\O_T} \varphi^2(x+t\theta)\left( p_{11}^{12}\mathcal{K}_1^+ \mathcal{K}_1^-
        +p_{12}^{12}\mathcal{K}_2^+ \mathcal{K}_1^-+p_{21}^{12}\mathcal{K}_1^+\mathcal{K}_2^-+p_{22}^{12}\mathcal{K}_2^+\mathcal{K}_2^-
        \right)~dxdt\\ + \int_{\O_T} \varphi(x+t\theta)\exp{(\mathrm{i}\lambda \phi)}\left( p_{11}^{12}\mathcal{K}_1^+ R^{v_0}_1
        +p_{12}^{12}\mathcal{K}_2^+ R^{v_0}_1+p_{21}^{12}\mathcal{K}_1^+R^{v_0}_2+p_{22}^{12}\mathcal{K}_2^+R^{v_0}_2
        \right)~dxdt\\ +\int_{\O_T} \varphi(x+t\theta)\exp{(-\mathrm{i}\lambda \phi)}\left( p_{11}^{12} \mathcal{K}_1^- R^{v_2}_1
        +p_{12}^{12} \mathcal{K}_2^- R^{v_2}_1+p_{21}^{12}\mathcal{K}_1^- R^{v_2}_2+p_{22}^{12}\mathcal{K}_2^- R^{v_2}_2
        \right)~dxdt\\+\int_{\O_T} \left( p_{11}^{12}R^{v_2}_1 R^{v_0}_1+p_{12}^{12}R^{v_2}_2 R^{v_0}_1+p_{21}^{12}R^{v_2}_1 R^{v_0}_2+p_{22}^{12}R^{v_2}_2 R^{v_0}_2
        \right)~dxdt=0,
    \end{aligned}
\end{align}
where $\overrightarrow{\mathcal{R}}^{v_2}_{\overrightarrow{\mathcal{K}}^{(+)}}=(R^{v_2}_1,R^{v_2}_2)^\top$ and $\overrightarrow{\mathcal{R}}^{v_0}_{\overrightarrow{\mathcal{K}}^{(-)}}=(R^{v_0}_1, R^{v_0}_2)^\top$ satisfying 
\begin{align}\label{both reminder term bounds}
\lambda\|\overrightarrow{\mathcal{R}}^{v_2}_{\overrightarrow{\mathcal{K}}^{(+)}}\|_{\textbf{L}^{2}(\Omega_{T})}+\|\partial_{t} \overrightarrow{\mathcal{R}}^{v_2}_{\overrightarrow{\mathcal{K}}^{(+)}}\|_{\textbf{L}^{2}(\Omega_{T})} &\leq C \|\varphi\|_{\textbf{H}^{3}{(\R^\ell})}\\\lambda\|\overrightarrow{\mathcal{R}}^{v_0}_{\overrightarrow{\mathcal{K}}^{(-)}}\|_{\textbf{L}^{2}(\Omega_{T})}+\|\partial_{t} \overrightarrow{\mathcal{R}}^{v_0}_{\overrightarrow{\mathcal{K}}^{(-)}}\|_{\textbf{L}^{2}(\Omega_{T})} &\leq C \|\varphi\|_{\textbf{H}^{3}{(\R^\ell})}.
\end{align}
Taking $\lambda\rightarrow\infty$ along with \eqref{both reminder term bounds}, we obtain that  
\begin{align} \label{Integral identity 1st linearization and after taking limit}
\int_{\O_T} \varphi^2(x+t\theta)\left( p_{11}^{12}\mathcal{K}_1^+ \mathcal{K}_1^-
        +p_{12}^{12}\mathcal{K}_2^+ \mathcal{K}_1^-+p_{21}^{12}\mathcal{K}_1^+\mathcal{K}_2^-+p_{22}^{12}\mathcal{K}_2^+\mathcal{K}_2^-
        \right)~dxdt=0,\end{align}  holds for any $\varphi \in \mathscr{C}_c^\infty(\mathbb{R}^\ell)$,  $\overrightarrow{\mathcal{K}}^{(+)}, \overrightarrow{\mathcal{K}}^{(-)}\in \mathbb{R}^2$  and $\theta\in \mathbb{S}^{\ell-1}$.
    Next to ensure that $p^{12}_{ij}=0$ in $\Omega$ for each $1\leq i,j\leq 2$,   we prescribe the specific choices of constant vectors $\overrightarrow{\mathcal{K}}^{(+)}$ and $\overrightarrow{\mathcal{K}}^{(-)}$. In particular, using \eqref{both reminder term bounds} along with  $\overrightarrow{\mathcal{K}}^{(+)}=\overrightarrow{\mathcal{K}}^{(-)}=(1,0)^\top$, \eqref{Integral identity 1st linearization and after taking limit} boils down to
\begin{align}\label{Intergal identity for p11}
    \begin{aligned}
        \int_{\O_T} \varphi^2(x+t\theta) p_{11}^{12}(x)~dxdt=0,~ \forall \varphi \in \mathscr{C}_c^\infty(\mathbb{R}^\ell) \text{ and } \theta\in \mathbb{S}^{\ell-1}.
    \end{aligned}
\end{align}
We extend $p_{11}^{12}$ by zero to $\mathbb{R}^\ell\setminus \O$, and continue to denote the resulting function by $p_{11}^{12}$. It follows that $p_{11}^{12}\in C_c^\infty(\mathbb{R}^\ell)$. Thus, we have
\begin{align}
      \int_{\mathbb{R}^\ell}\int_{0}^T \varphi^2(x+t\theta) p_{11}^{12}(x)~dxdt=0,~ \forall \varphi \in \mathscr{C}_c^\infty(\mathbb{R}^\ell) \text{ and } \theta\in \mathbb{S}^{\ell-1}.
\end{align}
The change  of variables $x+t\theta:=y$ then further leads us to
\begin{align}\label{Ray tranform where particular phi needed}
      \int_{\mathbb{R}^\ell}\int_{0}^T \varphi^2(y) p_{11}^{12}(y-t\theta)~dtdy=0,~ \forall \varphi \in \mathscr{C}_c^\infty(\mathbb{R}^\ell) \text{ and } \theta\in \mathbb{S}^{\ell-1}.
\end{align}
The next step is to invert the integral transform above by choosing a suitable $\varphi \in \mathscr{C}_c^\infty(\mathbb{R}^\ell)$. Although this step is standard and has been discussed in \cite{bhardwaj2026reconstructionpotentialdampingcoefficients, Isakov1991AnIH, kumar2026partial}, we present  here for the sake of completeness. 
Let $\rho \in C_{c}^{\infty}(\mathbb{B}_1(0))$ be a cutoff function satisfying $0\leq\rho\leq 1$ with $\int_{\mathbb{R}^\ell}\rho^2(x)dx=1$, where $\mathbb{B}_1(0)$ is a unit ball centered at $0$ and radius $1$. Also, let 
$\Omega_{r}\subset\mathbb{R}^\ell$ be a open set, given by 
\begin{align}\label{dom:Omega_epsilon}
    \Omega_{r}
    := \left\{ x \in \mathbb{R}^\ell \setminus \overline{\Omega} : 
    \operatorname{dist}(x, \Omega) < r \right\}, 
\end{align} 
for $r \in \left(0, \min\left\{1, \frac{T - \operatorname{diam}(\Omega)}{3}\right\}\right)$.
Further, for  $y_0\in \O_r$, we define \begin{align}\label{Phi_{delta}}
    \rho_{\delta}(y)=\frac{1}{\delta^{\frac{\ell}{2}}}\rho\left(\frac{y-y_0}{\delta}\right).
\end{align}
Then the sequence
$\{\rho_\delta\}_{\delta>0}\subset \mathscr{C}_{c}^{\infty}(\mathbb{R}^{\ell})$, whenever $0<\delta<r$. Using \eqref{Phi_{delta}} in \eqref{Ray tranform where particular phi needed} leads to
\begin{align}
    \frac{1}{\delta^\ell}\int_{B_\delta(y_0)}\int_0^T p_{11}^{12}(y-t\theta)\rho^2\left(\frac{y-y_0}{\delta}\right)\,dtdy=0, \text{ for all } 0<\delta<r \text{ and }y_0\in\O_r.
\end{align}
As $\delta\rightarrow\infty$, the aforementioned expression reduces to
\begin{align}\label{integral of a in term of T}
 \int_0^T p_{11}^{12}(y_0-t\theta)\,dt=0, \text{ for all }y_0\in\O_r  \text{ and }\theta\in\mathbb{S}^{\ell-1}. 
\end{align}
Further, we can be express the above integral as
\begin{align}\label{estimate a_{omega_{eps}}}
 \int_0^\infty p_{11}^{12}(y_0-t\theta)\,dt=0, \quad \text{for all }y_0\in\O_r \text{ and }\theta\in\mathbb{S}^{\ell-1},
\end{align}
because $p_{11}^{12}(y_0-t\theta)=0$, for every $t\geq T$. Indeed, because of $T>\diam(\O)$, one has 
\[
(\O_r \pm s\theta)\cap \O=\emptyset,\ \mbox{for all $s\geq T$ and $\theta\in\mathbb{S}^{\ell-1}$}.
\]
 From here, we aim to show that the $X-$ray transform of $p_{11}^{12}$ vanishes. To this end, we first observe that \eqref{estimate a_{omega_{eps}}} holds for all $\theta\in \mathbb{S}^{\ell-1}$, hence we arrive at 
\begin{align}\label{estimate 1 p11 on R}
 \int_{-\infty}^\infty p_{11}^{12}(y_0-t\theta)\,dt=0, \quad \text{for all }y_0\in\O_r \text{ and }\theta\in\mathbb{S}^{\ell-1}.
\end{align}
This tells us that the $X-$ray transform is vanishing on every line passing through $\O_r$. Following \cite{Rakesh01011988}, for any line $y+\mathbb{R}\theta$ where $\theta\in \mathbb{S}^{\ell-1}$ and $y\in \theta^\perp$, if $(y+\mathbb{R}\theta) \cap \O \neq\emptyset$, there exist $y_0 \in(y+\mathbb{R}\theta)\cap \O_r$, which results into 
\begin{align}\label{x ray}
 \int_{-\infty}^\infty p_{11}^{12}(y-t\theta)\,dt=0, \quad \text{for all }(\theta,y)\in T\mathbb{S}^{\ell-1}.
\end{align}
However, if $(y+\mathbb{R}\theta) \cap \O =\emptyset$, then \eqref{x ray} is trivially true, as $\supp (p_{11}^{12}) \subset \O$.
Finally, invoking the inversion of $X$-ray transform from (Theorem 1.1, \cite{Natterer86}) leads to the establishment of the unique recovery of the coefficients, that is $p_{11}^{12}=0$ in $\Omega$, which simply tells us that $p_{11}^1=p_{11}^2$ in $\Omega$. For the other coefficients of the matrix $\mathbb{P}^{12}$, again we need to choose the particular values of the constant vectors $\overrightarrow{\mathcal{K}}^{(+)}$ and $\overrightarrow{\mathcal{K}}^{(-)}$. Specifically, to recover $p_{ij}^{12}$, we set $\overrightarrow{\mathcal{K}}^{(+)} = e_j$ and $\overrightarrow{\mathcal{K}}^{(-)} = e_i$, where $e_i$ denotes the $i$-th standard basis vector. Subsequently, the preceding analysis completes the proof.  
\subsection{Unique determination of nonlinear potential}
\label{Recovery of the coefficients of nonlinearity} To do so, we first observe that the second-order partial derivative $\partial^2_{\epsilon_i \epsilon_j}$ does not lead to any extra information. So, to this end, 
we  apply $\partial_{\varepsilon_1 \varepsilon_2 \varepsilon_3}^3$ to the IBVP \eqref{govern eqn: For IBVP in epsilon form} and evaluate at $\varepsilon=0$, to get
\begin{equation}
\left\{\begin{array}{r l  c}
\displaystyle 
\mathcal{L}_{\mathbb{P}}\overrightarrow{w}(t,x) &= -\overrightarrow{b},\quad& (t,x)\in \Omega_T, \\[2pt]
\overrightarrow{w}(t,x)&= \overrightarrow{0},\quad&  (t,x)\in \partial\Omega_T ,\\[2pt]
\overrightarrow{w}(0,x)&= \partial_t \overrightarrow{w}(0,x)=\overrightarrow{0}, \quad& x\in \Omega,
\end{array}\right.
\end{equation}
where $\partial_{\varepsilon_1 \varepsilon_2 \varepsilon_3}^3 \overrightarrow{u}\rvert_{\varepsilon=0}=\overrightarrow{w}:=(w_1, w_2)^\top$, and
\begin{align}
\begin{aligned}
\overrightarrow{b}=
\begin{bmatrix}
   b_1\\[4pt]b_2
\end{bmatrix}=\begin{bmatrix}
  6q_{11}(x)v_1^{(1)}v_1^{(2)}v_1^{(3)}+2q_{12}(x)v_1^{(1)}v_1^{(3)}v_2^{(2)}+2q_{12}(x)v_1^{(1)}v_1^{(2)}v_2^{(3)}+2q_{12}(x)v_2^{(1)}v_1^{(2)}v_1^{(3)}\\[4pt]  2q_{21}(x)v_2^{(1)}v_2^{(3)}v_1^{(2)}+2q_{21}(x)v_2^{(1)}v_2^{(2)}v_1^{(3)}+2q_{21}(x)v_1^{(1)}v_2^{(2)}v_2^{(3)}+6q_{22}(x)v_2^{(1)}v_2^{(2)}v_2^{(3)}
\end{bmatrix}.
    \end{aligned}
\end{align}
A third-order linearization of the DN map \eqref{Def: DN map for IBVP of epsilon form} is
\begin{align}\label{Def: DN map 3-order linearization}
{}^3\Lambda_{\mathbb{P},\mathbb{Q}}\lvert_{\partial\O_T}:=\left(\partial^3_{\varepsilon_1 \varepsilon_2\varepsilon_3}\Lambda_{\mathbb{P},\mathbb{Q}}(\varepsilon\overrightarrow{f})\lvert_{\varepsilon=0}\right)\Big\lvert_{\partial\O_T}=\partial_\nu \overrightarrow{w}\lvert_{\partial\O_T}.
\end{align}
To prove uniqueness,  for $i=1,2$, we consider $\overrightarrow{w}^i$ be the solution to the IBVP
\begin{equation}
\left\{\begin{array}{r l  c}
\displaystyle \mathcal{L}_{\mathbb{P}}
\overrightarrow{w}^i(t,x) &=- \overrightarrow{b}^i,\quad& (t,x)\in \Omega_T, \\[2pt]
\overrightarrow{w}^i(t,x)&= \overrightarrow{0},\quad&  (t,x)\in \partial\Omega_T ,\\[2pt]
\overrightarrow{w}^i(0,x)&= \partial_t \overrightarrow{w}^i(0,x)= \overrightarrow{0}, \quad& x\in \Omega,
\end{array}\right.
\end{equation}
where 
\begin{align}
\begin{aligned}
\overrightarrow{b}^i=\begin{bmatrix}
    6q^i_{11}(x)v_1^{(1)}v_1^{(2)}v_1^{(3)}+2q^i_{12}(x)v_1^{(1)}v_1^{(3)}v_2^{(2)}+2q^i_{12}(x)v_1^{(1)}v_1^{(2)}v_2^{(3)}+2q_{12}^i(x)v_2^{(1)}v_1^{(2)}v_1^{(3)}\\[4pt]2q^i_{21}(x)v_2^{(1)}v_2^{(3)}v_1^{(2)}+2q^i_{21}(x)v_2^{(1)}v_2^{(2)}v_1^{(3)}+2q^i_{21}(x)v_1^{(1)}v_2^{(2)}v_2^{(3)}+6q^i_{22}(x)v_2^{(1)}v_2^{(2)}v_2^{(3)}
\end{bmatrix}.
    \end{aligned}
\end{align}
A third-order linearization of the DN map \eqref{Def: DN map for IBVP of epsilon form} corresponding to $\mathbb{P}=\mathbb{P}^i$ and $\mathbb{Q}=\mathbb{Q}^i$ is given by
\begin{align}\label{Def: DN map 3-order linearization with q=q1}
{}^3\Lambda_{\mathbb{P},\mathbb{Q}^i}\lvert_{\partial\O_T}=\partial_\nu \overrightarrow{w}^i\lvert_{\partial\O_T}.
\end{align}
The matrix $\mathbb{P}$ remains unchanged, since it's unique recovery has already been established. That's why a single notation $\mathbb{P} = \mathbb{P}^1 = \mathbb{P}^2$ is used. Also, we denote 
\begin{align}
    \overrightarrow{w}^{12}:=\overrightarrow{w}^1-\overrightarrow{w}^2,\quad \mathbb{Q}^{12}=\mathbb{Q}^1-\mathbb{Q}^2.
\end{align}
Then we have 
\begin{equation}\label{PDE: nonlinear simplified}
\left\{\begin{array}{r l  c}
\displaystyle \mathcal{L}_{\mathbb{P}}
 \overrightarrow{w}^{12}(t,x)&=- \overrightarrow{b}^{12},\quad& (t,x)\in \Omega_T,  \\[2pt]
\overrightarrow{w}^{12}(t,x)&= \overrightarrow{0},\quad&  (t,x)\in \partial\Omega_T ,\\[2pt]
\overrightarrow{w}^{12}(0,x)&= \partial_t \overrightarrow{w}^{12}(0,x)=\overrightarrow{0}, \quad& x\in \Omega,
\end{array}\right.
\end{equation}
where 
\begin{align}
\begin{aligned}
 \overrightarrow{b}^{12}=\begin{bmatrix}
    6q^{12}_{11}(x)v_1^{(1)}v_1^{(2)}v_1^{(3)}+2q^{12}_{12}(x)v_1^{(1)}v_1^{(3)}v_2^{(2)}+2q^{12}_{12}(x)v_1^{(1)}v_1^{(2)}v_2^{(3)}+2q_{12}^{12}(x)v_2^{(1)}v_1^{(2)}v_1^{(3)}\\[4pt]2q^{12}_{21}(x)v_2^{(1)}v_2^{(3)}v_1^{(2)}+2q^{12}_{21}(x)v_2^{(1)}v_2^{(2)}v_1^{(3)}+2q^{12}_{21}(x)v_1^{(1)}v_2^{(2)}v_2^{(3)}+6q^{12}_{22}(x)v_2^{(1)}v_2^{(2)}v_2^{(3)}
\end{bmatrix}.
    \end{aligned}
\end{align}
Using \(
      \Lambda_{\mathbb{P}^1,\mathbb{Q}^1}(\varepsilon\overrightarrow{f})= \Lambda_{\mathbb{P}^2,\mathbb{Q}^2}(\varepsilon\overrightarrow{f}),
\) in conjunction with \eqref{Def: DN map 3-order linearization with q=q1}, we have \(\partial_\nu \overrightarrow{w}^{12}\lvert_{\partial \O_T}=\overrightarrow{0}\).
Let $\overrightarrow{v_0}:=(v_{01},v_{02})^\top$ be the solution to the backward problem 
\begin{align}\label{PDE: backward problem used in integral identity nonlinear case}
    \begin{aligned}
      \left\{\begin{array}{l  c}
\displaystyle 
\Box \overrightarrow{v_0} +\mathbb{P}^\top\overrightarrow{v_0} = \overrightarrow{0},\quad& (t,x)\in \Omega_T,  \\[2pt]
\overrightarrow{v_0}(T,x)= \partial_t \overrightarrow{v_0}(T,x)=\overrightarrow{0}, \quad& x\in \Omega.
\end{array}\right.
    \end{aligned}
\end{align}
Multiply \eqref{PDE: nonlinear simplified} by $\overrightarrow{v_0}$ and integrate over $\O_T$ to get
\begin{align}
    &-\int_{\O_T} \left(\Box \overrightarrow{w}^{12} +\mathbb{P}\overrightarrow{w}^{12}\right)\cdot \overrightarrow{v_0}~dxdt \\&=     \int_{\O_T}(b^{12}_1 v_{01}+b_{2}^{12}v_{02})~dxdt\\&=   \int_{\O_T} \big(6q^{12}_{11}(x)v_1^{(1)}v_1^{(2)}v_1^{(3)}+2q^{12}_{12}(x)v_1^{(1)}v_1^{(3)}v_2^{(2)}+2q^{12}_{12}(x)v_1^{(1)}v_1^{(2)}v_2^{(3)}+2q_{12}^{12}(x)v_2^{(1)}v_1^{(2)}v_1^{(3)}\big) v_{01}dxdt\\& +\int_{\O_T} \big(2q^{12}_{21}(x)v_2^{(1)}v_2^{(3)}v_1^{(2)}+2q^{12}_{21}(x)v_2^{(1)}v_2^{(2)}v_1^{(3)}+2q^{12}_{21}(x)v_1^{(1)}v_2^{(2)}v_2^{(3)}+6q^{12}_{22}(x)v_2^{(1)}v_2^{(2)}v_2^{(3)}
    \big) v_{02}dxdt,
\end{align}
where $\overrightarrow{b}^{12}:=(b_1^{12}, b_2^{12})^\top.$
The integration by parts formula, along with initial/final time information and boundary data in \eqref{PDE: nonlinear simplified} and \eqref{PDE: backward problem used in integral identity nonlinear case} leads to
\begin{align}\label{Integral: component wise}
\begin{aligned}
    & \int_{\O_T} \big(6q^{12}_{11}(x)v_1^{(1)}v_1^{(2)}v_1^{(3)}+2q^{12}_{12}(x)v_1^{(1)}v_1^{(3)}v_2^{(2)}+2q^{12}_{12}(x)v_1^{(1)}v_1^{(2)}v_2^{(3)}+2q_{12}^{12}(x)v_2^{(1)}v_1^{(2)}v_1^{(3)}\big) v_{01}~dxdt\\&\qquad+\int_{\O_T} \big(2q^{12}_{21}(x)v_2^{(1)}v_2^{(3)}v_1^{(2)}+2q^{12}_{21}(x)v_2^{(1)}v_2^{(2)}v_1^{(3)}+2q^{12}_{21}(x)v_1^{(1)}v_2^{(2)}v_2^{(3)}\\& \qquad\qquad\qquad+6q^{12}_{22}(x)v_2^{(1)}v_2^{(2)}v_2^{(3)}
    \big) v_{02}~dxdt=0.
    \end{aligned}
\end{align} Next, we choose (see Proposition \ref{Asymptotic solutions for IBVP} and \ref{Asymptotic solutions for backward IBVP})
\begin{equation}\label{eqn: IVP geom opt sol 1 in nonlinear id}
 \overrightarrow{v}^{(i)}(t,x)= \exp{(\mathrm{i}\lambda \phi)}\bigg(z_{0}(t,x)\overrightarrow{K}^{+,i}+\frac{Z_{1,\overrightarrow{K}^{+,i}}(t,x)}{\lambda}+\frac{Z_{2,\overrightarrow{K}^{+,i}}(t,x)}{\lambda^2}+\cdots+\frac{Z_{N,\overrightarrow{K}^{+,i}}(t,x)}{\lambda^N}\bigg)+\overrightarrow{\mathcal{R}}_{\overrightarrow{K}^{+,i}},
    \end{equation}
 and 
  \begin{equation}\label{eqn: ansatz2 nonlinear id}
 \overrightarrow{v_0}(t,x)= \exp{(-3\mathrm{i}\lambda \phi)}\bigg(w_{0}(t,x)\overrightarrow{K}^{-}+\frac{W_{1,\overrightarrow{K}^{-}}(t,x)}{\lambda}+\frac{W_{2,\overrightarrow{K}^{-}}(t,x)}{\lambda^2}+\cdots+\frac{W_{N,\overrightarrow{K}^{-}}(t,x)}{\lambda^N}\bigg)+\overrightarrow{\mathcal{R}}_{\overrightarrow{K}^{-}},
    \end{equation}
 where $\overrightarrow{v}^{(i)}$ and $\overrightarrow{v_0}$ are the asymptotic solutions to \eqref{PDE: v power i} and \eqref{PDE: backward problem used in integral identity nonlinear case}, respectively. Let $\varphi \in \mathscr{C}_c^\infty(\mathbb{R}^\ell)$ then $z_0=w_0=\varphi(x+t\theta)$ solve the transport equation \eqref{IVP: transport equation z_0} and \eqref{IVP: transport w_0}. The remainder terms  $\overrightarrow{\mathcal{R}}_{\overrightarrow{K}^{+,i}}(t,x,\lambda)$ for $1\leq i\leq 3$ and $\overrightarrow{\mathcal{R}}_{\overrightarrow{K}^{-}}(t,x,\lambda)$, satisfy the following estimates: \begin{align}\label{Bound: correction terms}\lVert \overrightarrow{\mathcal{R}}_{\overrightarrow{K}^{+,i}}\rVert_{\textbf{L}^{\infty}(\O_T)}\leq C\lambda^{-1} \quad \mbox{ and}\quad \lVert \overrightarrow{\mathcal{R}}_{\overrightarrow{K}^{-}}\rVert_{\textbf{L}^{\infty}(\O_T)}\leq C\lambda^{-1},\end{align} 
   where $C>0$ is independent of $\lambda$. Also, we have 
    \begin{align}  \overrightarrow{\mathcal{R}}_{\overrightarrow{K}^{+,i}}(0,x,\lambda)=\partial_t \overrightarrow{\mathcal{R}}_{\overrightarrow{K}^{+,i}}(0,x,\lambda)=\overrightarrow{0}  \mbox{ and }\; \overrightarrow{\mathcal{R}}_{\overrightarrow{K}^{-}}(T,x,\lambda)=\partial_{t}\overrightarrow{\mathcal{R}}_{\overrightarrow{K}^{-}}(T,x,\lambda)=\overrightarrow{0},\ \end{align}
    for $x\in \O$ and $1\leq i\leq 3$. For $1\leq j\leq N$ and $1\leq i\leq 3$, we
    denote 
    \begin{align}
        \begin{aligned}
            z_0 \overrightarrow{K}^{+,i}=\begin{bmatrix}
                z_0 K^{+,i}_1\\ z_0 K^{+,i}_2 
            \end{bmatrix},~~  w_0 \overrightarrow{K}^{+,i}=\begin{bmatrix}
                w_0 K^{-}_1\\ w_0 K^{-}_2 
            \end{bmatrix}, ~~Z_{j,\overrightarrow{K}^{+,i}}=\begin{bmatrix}
Z_{1,(j,\overrightarrow{K}^{+,i})}\\Z_{2,(j,\overrightarrow{K}^{+,i})}  
            \end{bmatrix},
   \\[4pt]
          W_{j,\overrightarrow{K}^-}=\begin{bmatrix}
W_{1,(j,\overrightarrow{K}^-)}\\W_{2,(j,\overrightarrow{K}^-)}
            \end{bmatrix}, ~~ \overrightarrow{\mathcal{R}}_{\overrightarrow{K}^{+,i}}=\begin{bmatrix}
R_{1,\overrightarrow{K}^{+,i}}\\R_{2,\overrightarrow{K}^{+,i}}
            \end{bmatrix}, ~~\overrightarrow{\mathcal{R}}_{\overrightarrow{K}^{-}}=\begin{bmatrix}
R_{1,\overrightarrow{K}^{-}}\\R_{2,\overrightarrow{K}^{-}}
            \end{bmatrix}, \\[4pt]\mathcal{Z}_i=\begin{bmatrix}
             \mathcal{Z}_{1,i}\\[2pt]
\mathcal{Z}_{2,i}\end{bmatrix}:=\begin{bmatrix}
          \displaystyle  \sum_{j=1}^N \dfrac{Z_{1,(j,\overrightarrow{K}^{+,i})}}{\lambda^j}\\[4pt]\displaystyle\sum_{j=1}^N \dfrac{Z_{2,(j,\overrightarrow{K}^{+,i})}}{\lambda^j}
            \end{bmatrix}, \text{ and }\mathcal{W}=\begin{bmatrix}
\mathcal{W}_1\\\mathcal{W}_2  
\end{bmatrix}:=\begin{bmatrix}
               \displaystyle \sum_{j=1}^N \dfrac{W_{1,(j,\overrightarrow{K}^{-})}}{\lambda^j}\\[2pt]\displaystyle\sum_{j=1}^N \dfrac{W_{2,(j,\overrightarrow{K}^{-})}}{\lambda^j}
            \end{bmatrix}.
        \end{aligned}
    \end{align}  
   For $1\leq k,l\leq 2$, the above notation in conjunction with \eqref{Bound: correction terms} leads to \\[2pt]
    \textbf{Case~1:} $k\neq l$\\
    \begin{align}
        \begin{aligned}
v_k^{(1)}v_k^{(2)}v_l^{(3)}v_{0k}=\varphi^4(x+t\theta)K^{+,1}_k K^{+,2}_k K^{+,3}_l K^{-}_k+\mathcal{O}(\lambda^{-1}),\\[2pt]v_k^{(1)}v_l^{(2)}v_k^{(3)}v_{0k}=\varphi^4(x+t\theta)K^{+,1}_k K^{+,2}_l K^{+,3}_k K^{-}_k+\mathcal{O}(\lambda^{-1}),\\v_l^{(1)}v_k^{(2)}v_k^{(3)}v_{0k}=\varphi^4(x+t\theta)K^{+,1}_l K^{+,2}_k K^{+,3}_k K^{-}_k+\mathcal{O}(\lambda^{-1}),
        \end{aligned}
    \end{align}
     \textbf{Case~2:} $k= l$\\
    \begin{align}
        \begin{aligned}
v_k^{(1)}v_k^{(2)}v_k^{(3)}v_{0k}=\varphi^4(x+t\theta)K^{+,1}_k K^{+,2}_k K^{+,3}_k K^{-}_k+\mathcal{O}(\lambda^{-1}).
        \end{aligned}
    \end{align}
Here we have used that $z_0=w_0= \varphi(x+t\theta)$, where $\varphi \in \mathscr{C}_c^\infty(\mathbb{R}^\ell)$. Order $\lambda^{-1}$ appears from the residual terms beyond the leading term in the aforementioned asymptotic expansion and/or from the estimates correction terms, possibly in combination with other lower-order terms.
The substitution of aforementioned GO solution in \eqref{Integral: component wise} and use of the above expression yields
\begin{align}\label{Integral identity: nonlinear}
    \begin{aligned}
\int_{\O_T}\varphi^4(x+t\theta)&\big(6q^{12}_{11}(x)K^{+,1}_1 K^{+,2}_1 K^{+,3}_1 K^{-}_1+2q^{12}_{12}(x)K^{+,1}_1 K^{+,2}_2 K^{+,3}_1 K^{-}_1\\[2pt]&+2q^{12}_{12}(x)K^{+,1}_1 K^{+,2}_1 K^{+,3}_2 K^{-}_1+2q^{12}_{12}(x) K^{+,1}_2 K^{+,2}_1 K^{+,3}_1 K^{-}_1\\[2pt]&+2q_{21}^{12}(x)K^{+,1}_2  K^{+,2}_1 K^{+,3}_2 K^{-}_2+2q_{21}^{12}(x)K^{+,1}_2 K^{+,2}_2 K^{+,3}_1 K^{-}_2\\[2pt]&+2q_{21}^{12}(x)K^{+,1}_1 K^{+,2}_2 K^{+,3}_2 K^{-}_2+6q_{22}^{12}(x)K^{+,1}_2 K^{+,2}_2 K^{+,3}_2 K^{-}_2
\big)
        dxdt=0,
    \end{aligned}
\end{align}
as $\lambda\rightarrow \infty.$ To ensure the unique identification of each coefficient of nonlinearity, we prescribe a specific choice of the constant vectors $\overrightarrow{K}^{+,i}$ and $\overrightarrow{K}^{-}$. In particular, we choose $\overrightarrow{K}^{+,1}=\overrightarrow{K}^{+,2}=\overrightarrow{K}^{+,3}=\overrightarrow{K}^{-}=e_1$, where $e_i$ denotes the $i$-th standard basis vector of $\mathbb{R}^2$. 
Using these choices of constant vector in \eqref{Integral identity: nonlinear} leads to
\begin{align}
    \int_{\O_T}\varphi^4(x+t\theta)q^{12}_{11}(x)~dxdt=0,~ \forall \varphi \in \mathscr{C}_c^\infty(\mathbb{R}^\ell) \text{ and } \theta\in \mathbb{S}^{\ell-1}.
\end{align}
We implement an argument analogous to that used to establish the uniqueness of the coefficient from the integral identity \eqref{Intergal identity for p11} to show that $q_{11}^{12}=0$. However, establishing the unique recovery of the coefficients of the nonlinearity requires a suitable choice of $\varphi \in \mathscr{C}_c^\infty(\mathbb{R}^\ell)$. Extend $q_{11}^{12}$ by zero to $\mathbb{R}^\ell\setminus \O$, and use the change of  variable to get
\begin{align}\label{Ray tranform where particular phi needed but in nonlinear}
      \int_{\mathbb{R}^\ell}\int_{0}^T \varphi^4(y) q_{11}^{12}(y-t\theta)~dtdy=0,~ \forall \varphi \in \mathscr{C}_c^\infty(\mathbb{R}^\ell) \text{ and } \theta\in \mathbb{S}^{\ell-1}.
\end{align}
For $y_0\in \O_r$, we define \begin{align}\label{Phi_{delta}1}
    \rho_{\delta}(y)=\frac{1}{\delta^{\frac{\ell}{4}}}\rho\left(\frac{y-y_0}{\delta}\right),
\end{align}
where $\rho \in C_{c}^{\infty}(\mathbb{B}_1(0))$ is a cutoff function lies in $[0,1]$ with $\int_{\mathbb{R}^\ell}\rho^4(x)dx=1$.
It follows from the above construction that the sequence
$\{\rho_\delta\}_{\delta>0}\subset \mathscr{C}_{c}^{\infty}(\mathbb{R}^{\ell})$ whenever $0<\delta<r$.  Thus, we have
\begin{align}
    \frac{1}{\delta^\ell}\int_{B_\delta(y_0)}\int_0^T q_{11}^{12}(y-t\theta)\rho^4\left(\frac{y-y_0}{\delta}\right)\,dtdy=0, \text{ for all } 0<\delta<r \text{ and }y_0\in\O_r,
\end{align}
which then leads to 
\begin{align}\label{integral of q11 in term of T}
 \int_0^T q_{11}^{12}(y_0-t\theta)\,dt=0, \text{ for all }y_0\in\O_r  \text{ and }\theta\in\mathbb{S}^{\ell-1}, 
\end{align}
as $\delta\rightarrow\infty$. The foregoing integral is now reduced to \eqref{integral of a in term of T} from which we concluded that $p_{11}^{12}=0$, the only change is that the coefficient $q_{11}^{12}$ appears in place of $p_{11}^{12}$. Since both coefficients belong to the same function class $C_c^\infty(\O)$, we  have $q_{11}^{12}=0$ in $\Omega$, which simply tells that $q_{11}^1=q_{11}^2$ in $\Omega$. For the other coefficients of the matrix $\mathbb{Q}^{12}$, again we need to choose the particular values of the constant vectors $\overrightarrow{K}^{+,i}$ and $\overrightarrow{K}^{-}$. Specifically, to recover the remaining coefficients, we consider the suitable choices of constant vectors given in Table \ref{table1}. 
\begin{table}[ht]
\centering
\begin{tabular}{|Sc|Sc|Sc|Sc|Sc|}
   \hline ~ & $\overrightarrow{K}^{+,1}$&$\overrightarrow{K}^{+,2}$&$\overrightarrow{K}^{+,3}$&$\overrightarrow{K}^{-}$\\\hline
     $q_{12}^{12}$& $e_1$&$e_2$&$e_1$&$e_1$\\\hline $q_{21}^{12}$& $e_2$&$e_1$&$e_2$&$e_2$\\\hline$q_{22}^{12}$& $e_2$&$e_2$&$e_2$&$e_2$\\\hline
\end{tabular}
\vspace{0.5cm}
\caption{A particular choice of $\overrightarrow{K}^{+,i}$ and $\overrightarrow{K}^{-}$}
\label{table1}
\end{table}
Subsequently, the preceding analysis completes the proof.   
\section{Proof of Theorem \ref{thm:main_result_partial_data}}\label{Recovery of the coefficients from partial DN map} 
In the following, we establish the unique identification of the coefficients appearing in \eqref{govern eqn: For IBVP in epsilon form} from the partial DN map given by 
    \begin{align}
\Lambda^{\sharp}_{\mathbb{P},\mathbb{Q}}(\varepsilon\overrightarrow{f}):=\partial_\nu \overrightarrow{u}\rvert_{(0,T)\times\Gamma},
    \end{align}
    where $\overrightarrow{u}$ is a solution to \eqref{govern eqn: For IBVP in epsilon form} and $\Gamma\subset \partial \Omega$ is a nonempty open set.  First, in \S \ref{Subsection: Unique continuation principle}, we recall the unique continuation principle established in \cite{BenFraj2025}, which is then used in converting the partial-data problem to a full-data problem.
    Next, in \S \ref{Recovery of matrix-valued potential from partial DN map} and \S \ref{Recovery of nonlinear potential from partial DN map}, we discuss the unique identification of matrix-valued potential and nonlinear potential through the partial DN map, respectively.

\subsection{Unique continuation principle}\label{Subsection: Unique continuation principle}
We begin with defining some notations to state the unique continuation principle (UCP) associated with the linear hyperbolic system \begin{equation}\label{PDE: UCP}
\left\{\begin{array}{r l c}
\displaystyle \mathcal{L}_{\mathbb{P}}\overrightarrow{u}(t,x)& = \overrightarrow{F},& (t,x)\in \Omega_T,  \\[2pt]
\overrightarrow{u}(t,x)&= \overrightarrow{0},&  (t,x)\in \partial\Omega_T ,\\[2pt]
\overrightarrow{u}(0,x)&= \overrightarrow{0}, \quad\partial_t\overrightarrow{u}(0,x)= \overrightarrow{0}, &\quad x\in \Omega,
\end{array}\right.
\end{equation}
where $\mathbb{P}\in \mathcal{P}$ and $\overrightarrow{F} \in \textbf{H}^2(\O_T)$. 
For the  open set $\mathcal{U}$ in $\overline{\Omega}$, we introduce a family of relatively compact open subsets of $\mathcal{U}$ in $\overline{\Omega}$. To be more precise, we consider the  open sets $\mathcal{U}_{j} \subseteq \mathcal{U}$, $j=1,2,3$, such that $\overline{\mathcal{U}}_{j+1} \subset \mathcal{U}_{j}$ for $j=1,2$ and  $\partial\Omega\subset\partial\mathcal{U}_{j}$ for $j=1,2,3$.   
We also assume that
$\overline{\mathcal U}_1\cap(\partial\mathcal U\setminus\partial\Omega)
=\emptyset$.
  Now, we recall the following results from \cite{BenFraj2025}.
\begin{theorem}\label{UCP} \cite[Lemma~~3.1]{BenFraj2025}
    Let $\Gamma\subset \partial \Omega$ be a nonempty open set, $T>0$ be sufficiently large and $\mathbb{P}\in \mathcal{P}$. Also, let $\overrightarrow{u}$ be the solution to the IBVP \eqref{PDE: UCP}. Then there exist positive constants $T^{'}<T$, $C$, $A$ and $\eta_0$ such that 
    \begin{align}
        \lVert \overrightarrow{u}\rVert_{\textbf{H}^1((0,T^{'})\times(\mathcal{U}_2\setminus\mathcal{U}_3))} \leq \frac{C}{\eta^{1/2}}\lVert \overrightarrow{u}\rVert_{\textbf{H}^2(\O_T)}+e^{A\eta}\left(\lVert \overrightarrow{F}\rVert_{\textbf{L}^2{((0,T)\times \mathcal{U})}}+\lVert \partial_\nu \overrightarrow{u}\rVert_{\textbf{L}^2((0,T)\times \Gamma)}
        \right),
    \end{align}
    holds for any $\eta>\eta_0$.
\end{theorem}
\begin{corollary}\label{Corollary}
    In addition to the assumptions of Theorem \ref{UCP}, suppose that \(\operatorname{supp}(\overrightarrow{F}) \subset (0,T)\times (\Omega\setminus \mathcal{U})\) and \(\partial_{\nu}\overrightarrow{u}=\overrightarrow{0}\) on \((0,T)\times\Gamma\). Then \(\overrightarrow{u}=\overrightarrow{0}\) in \((0,T^{'})\times (\mathcal{U}_2\setminus \mathcal{U}_3)\).
\end{corollary}
    \subsection{Recovery of matrix-valued potential}\label{Recovery of matrix-valued potential from partial DN map}
   In this subsection, we give a proof for converting  the partial data problem to a full data problem.  We start by recalling the IBVP (see \eqref{PDE: used in integral identity}) \begin{align}\label{PDE: used in integral identity nonlinear}
    \left\{\begin{array}{l  c}
\displaystyle 
\Box \overrightarrow{v}^{12} +\mathbb{P}^1\overrightarrow{v}^{12} = -\mathbb{P}^{12}\overrightarrow{v}^2,\quad& (t,x)\in \Omega_T,  \\[2pt]
\overrightarrow{v}^{12}(t,x)= \overrightarrow{0},\quad&  (t,x)\in \partial\Omega_T ,\\[2pt]
\overrightarrow{v}^{12}(0,x)= \partial_t \overrightarrow{v}^{12}(0,x)=\overrightarrow{0}, \quad& x\in \Omega,
\end{array}\right.
\end{align}
where $\overrightarrow{v}^j$ is solution to \eqref{PDE: compact form vi}, $\overrightarrow{v}^{12}=\overrightarrow{v}^1-\overrightarrow{v}^2$, and $\mathbb{P}^{12}=\mathbb{P}^1-\mathbb{P}^2$.
Since \(
      \Lambda_{\mathbb{P}^1,\mathbb{Q}^1}^\sharp(\varepsilon\overrightarrow{f})= \Lambda_{\mathbb{P}^2,\mathbb{Q}^2}^\sharp(\varepsilon\overrightarrow{f}),
\) the first order linearization of partial DN map leads to
\[{}^1\Lambda_{\mathbb{P}^1}^\sharp(f,g)\lvert_{(0,T)\times \Gamma}={}^1\Lambda_{\mathbb{P}^2}^\sharp(f,g)\lvert_{(0,T)\times \Gamma}, \quad f,g\in H^{m+1}(\partial\O_T).\]
Consequently, we have $\partial_\nu \overrightarrow{v}^{12}\lvert_{(0,T)\times \Gamma}=\overrightarrow{0}$.
Let $\Xi \in C^{\infty}(\overline\Omega,[0,1])$ be the smooth cutoff function given by
\begin{equation}\label{eqn: construction theta}
\Xi := \begin{cases}
0 \qquad\text{ in } \mathcal{U}_{3},\\
1 \qquad\text{ in } \Omega\setminus\overline{\mathcal{U}_2}.
\end{cases}
\end{equation}
Set $\overrightarrow{\widetilde{v}}:=\Xi  \overrightarrow{v}^{12}$, then $\overrightarrow{\widetilde{v}}$ solves
\begin{align}\label{v tilde equation}
    \begin{cases}
      \Box \overrightarrow{\widetilde{v}} +\mathbb{P}^1 \overrightarrow{\widetilde{v}} = \Xi \left(-\mathbb{P}^{12}\overrightarrow{v}^2\right)- (\Delta \Xi ) \overrightarrow{\widetilde{v}}-2 \nabla  \Xi \cdot\nabla \overrightarrow{\widetilde{v}}  ,  & (t,x) \in \O_T,\\
      \overrightarrow{\widetilde{v}}(t,x)  = 0, & (t,x) \in\partial \O_T,\\
      \overrightarrow{\widetilde{v}}(0, x )=0, \quad \partial_t \overrightarrow{\widetilde{v}}(0, x ) = 0,  &\quad  x \in\Omega,
    \end{cases}
\end{align}
where $\overrightarrow{v}^{12}$ is a solution to IBVP \eqref{PDE: used in integral identity nonlinear}.
Since $\mathbb{P}^1=\mathbb{P}^2$ in $\mathcal{U}$, we have $\mathbb{P}^{12}=0$ in $\mathcal{U}$. From the construction given by \eqref{eqn: construction theta}, we have $\Xi  \mathbb{P}^{12}=\mathbb{P}^{12}$ in $\overline{\Omega}$. In addition, we have $- (\Delta \Xi ) \overrightarrow{\widetilde{v}}-2 \nabla  \Xi \cdot\nabla \overrightarrow{\widetilde{v}}\equiv \overrightarrow{0}$ in $(0,T)\times (\mathcal{U}_3\cup (\Omega\setminus\overline{\mathcal{U}_2}))$. Thus, $\overrightarrow{F}:= \Xi \left(-\mathbb{P}^{12}\overrightarrow{v}^2\right)- (\Delta \Xi ) \overrightarrow{\widetilde{v}}-2 \nabla  \Xi \cdot\nabla \overrightarrow{\widetilde{v}}$ is supported in $(0,T)\times (\O\setminus \mathcal{U})$. 
Since $\overrightarrow{\widetilde{v}}\equiv\overrightarrow{0}$ in $\mathcal{U}_3$, we have 
\(
\partial_\nu \overrightarrow{\widetilde{v}}|_{\partial\O_T}=0.
\) From Corollary \ref{Corollary}, we then have $\overrightarrow{\widetilde{v}}=\overrightarrow{0}\text{ in }((0,T^{'})\times (\mathcal{U}_2\setminus \mathcal{U}_3)$. Next, we multiply \eqref{v tilde equation} with $\overrightarrow{v_0}$, where $\overrightarrow{v_0}$ solves 
\begin{align}
\begin{cases}
    \Box \overrightarrow{v_0}+(\mathbb{P}^1)^\top \overrightarrow{v_0}=\overrightarrow{0},& (t,x) \in \O_{T^{'}}:=(0,T')\times\O,\\ 
    \overrightarrow{v_0}(T^{'},x)=\partial_t \overrightarrow{v_0}(T^{'},x)=\overrightarrow{0}, & \qquad \qquad x \in \O,
    \end{cases}
\end{align}
and integrate over $\O_{T^{'}}$ to get
\begin{align}
     \int_{\O_{T^{'}}} \left(
    \mathbb{P}^{12}\overrightarrow{v}^2
    \right)\cdot \overrightarrow{v_0}~dxdt=\overrightarrow{0}.
\end{align}
Choose (see Lemma \ref{lemma: geom opt sol v_{1}} and Lemma \ref{lemma: geom opt sol v})
\begin{align}\label{GO solution v2 non}
    \overrightarrow{v}^{2}(t,x) = \varphi(x+t\theta)\exp{(\mathrm{i}\lambda \phi)}\overrightarrow{\mathcal{K}}^{(+)} + \overrightarrow{\mathcal{R}}^{v_2}_{\overrightarrow{\mathcal{K}}^{(+)}}(t,x,\lambda),
\end{align}
and 
\begin{align}\label{GO solution y non}
     \overrightarrow{v_0}(t,x) = \varphi(x+t\theta)\exp{(-\mathrm{i}\lambda \phi)}\overrightarrow{\mathcal{K}}^{(-)} + \overrightarrow{\mathcal{R}}^{v_0}_{\overrightarrow{\mathcal{K}}^{(-)}}(t,x,\lambda),
\end{align}
for all $\varphi\in \mathtt{C}_c^{\infty}(\mathbb{R}^\ell)$, where 
\begin{align}\label{eq; definition of C(Rl)q}
\mathtt{C}_c^{\infty}(\mathbb{R}^\ell):=\{\varphi\in C_c^{\infty}(\mathbb{R}^\ell):\supp(\varphi)\cap\Omega=\emptyset\ \mbox{and}\ (\supp(\varphi)\pm T^{'}\theta)\cap \Omega=\emptyset,\ \forall \ \theta\in\mathbb{S}^{\ell-1}\}.
\end{align}
Here the remainder term $\overrightarrow{\mathcal{R}}_{\overrightarrow{\mathcal{K}}^{(+)}}^{v_2}(t,x,\lambda)$ satisfies
\begin{align}\label{R_{1} initial}
   \overrightarrow{\mathcal{R}}_{\overrightarrow{\mathcal{K}}^{(+)}}^{v_2}(0,x,\lambda) = \partial_{t}\overrightarrow{\mathcal{R}}_{\overrightarrow{\mathcal{K}}^{(+)}}^{v_2}(0,x,\lambda) =  \overrightarrow{0}   \text{ in } \Omega,
\end{align}
with 
\begin{align}
\lambda\|\overrightarrow{\mathcal{R}}_{\overrightarrow{\mathcal{K}}^{(+)}}^{v_2}\|_{\textbf{L}^{2}(\Omega_{T})}+\|\partial_{t} \overrightarrow{\mathcal{R}}_{\overrightarrow{\mathcal{K}}^{(+)}}^{v_2}\|_{\textbf{L}^{2}(\Omega_{T})} \leq C \|\varphi\|_{\textbf{H}^{3}{(\R^\ell})}.
\end{align}
From the last relation, we observe that 
\begin{align}
\lambda\|\overrightarrow{\mathcal{R}}_{\overrightarrow{\mathcal{K}}^{(+)}}^{v_2}\|_{\textbf{L}^{2}(\Omega_{T^{'}})}+\|\partial_{t} \overrightarrow{\mathcal{R}}_{\overrightarrow{\mathcal{K}}^{(+)}}^{v_2}\|_{\textbf{L}^{2}(\Omega_{T^{'}})} \leq C \|\varphi\|_{\textbf{H}^{3}{(\R^\ell})}.
\end{align}
Also, the correction term $\overrightarrow{\mathcal{R}}^{v_0}_{\overrightarrow{\mathcal{K}}^{(-)}}$ satisfies 
\begin{align}\label{R_{1} final times}
\overrightarrow{\mathcal{R}}^{v_0}_{\overrightarrow{\mathcal{K}}^{(-)}}(T',x,\lambda) = \partial_{t}\overrightarrow{\mathcal{R}}^{v_0}_{\overrightarrow{\mathcal{K}}^{(-)}}(T',x,\lambda) =  \overrightarrow{0}   \text{ in } \Omega,
\end{align}
with
\begin{align}\label{eq: estimate of r_{(2)} non}
\lambda\|\overrightarrow{\mathcal{R}}^{v_0}_{\overrightarrow{\mathcal{K}}^{(-)}}\|_{\textbf{L}^{2}(\Omega_{T^{'}})}+\|\partial_{t} \overrightarrow{\mathcal{R}}^{v_0}_{\overrightarrow{\mathcal{K}}^{(-)}}\|_{\textbf{L}^{2}(\Omega_{T'})} \leq C \|\varphi\|_{\textbf{H}^{3}{(\R^\ell})}.
\end{align}
Proceeding as in \S \ref{Recovery of zeroth-order coefficients}, we obtain the unique identification of $\mathbb{P}$, that is, $\mathbb{P}^1 = \mathbb{P}^2$ in $\O$.
 \subsection{Recovery of nonlinear potential}\label{Recovery of nonlinear potential from partial DN map}
 We start by recalling the IBVP (see \eqref{PDE: nonlinear simplified})
\begin{equation}\label{PDE: nonlinear simplified partial}
\left\{\begin{array}{r l  c}
\displaystyle \mathcal{L}_{\mathbb{P}} \overrightarrow{w}^{12}(t,x)&=- \overrightarrow{b}^{12},\quad& (t,x)\in \Omega_T,  \\[2pt]
\overrightarrow{w}^{12}(t,x)&= \overrightarrow{0},\quad&  (t,x)\in \partial\Omega_T ,\\[2pt]
\overrightarrow{w}^{12}(0,x)&= \partial_t \overrightarrow{w}^{12}(0,x)=\overrightarrow{0}, \quad& x\in \Omega,
\end{array}\right.
\end{equation}
where 
\begin{align}
\begin{aligned}
 \overrightarrow{b}^{12}=\begin{bmatrix}
    6q^{12}_{11}(x)v_1^{(1)}v_1^{(2)}v_1^{(3)}+2q^{12}_{12}(x)v_1^{(1)}v_1^{(3)}v_2^{(2)}+2q^{12}_{12}(x)v_1^{(1)}v_1^{(2)}v_2^{(3)}+2q_{12}^{12}(x)v_2^{(1)}v_1^{(2)}v_1^{(3)}\\[4pt]2q^{12}_{21}(x)v_2^{(1)}v_2^{(3)}v_1^{(2)}+2q^{12}_{21}(x)v_2^{(1)}v_2^{(2)}v_1^{(3)}+2q^{12}_{21}(x)v_1^{(1)}v_2^{(2)}v_2^{(3)}+6q^{12}_{22}(x)v_2^{(1)}v_2^{(2)}v_2^{(3)}
\end{bmatrix},
    \end{aligned}
\end{align}
$\mathbb{P} = \mathbb{P}^1 = \mathbb{P}^2$, 
\( \overrightarrow{w}^{12}=\overrightarrow{w}^1-\overrightarrow{w}^2,\text{ and } \mathbb{Q}^{12}=\mathbb{Q}^1-\mathbb{Q}^2.\)
Since \(
      \Lambda_{\mathbb{P}^1,\mathbb{Q}^1}^\sharp(\varepsilon\overrightarrow{f})= \Lambda_{\mathbb{P}^2,\mathbb{Q}^2}^\sharp(\varepsilon\overrightarrow{f})
\), the third-order linearization of the partial DN map leads to $\partial_\nu \overrightarrow{w}^{12}|_{(0,T)\times \Gamma}=\overrightarrow{0}$. 
Set $\overrightarrow{\widetilde{w}}:=\Xi  \overrightarrow{w}^{12}$, then $\overrightarrow{\widetilde{w}}$ solves the IBVP
\begin{align}\label{PDE: tilde w}
 \left\{\begin{array}{l  c}
\displaystyle \Box \overrightarrow{\widetilde{w}} +\mathbb{P} \overrightarrow{\widetilde{w}}=- \Xi \overrightarrow{b}^{12}- (\Delta \Xi ) \overrightarrow{\widetilde{w}}-2 \nabla  \Xi \cdot\nabla \overrightarrow{\widetilde{w}} ,\quad& (t,x)\in \Omega_T,  \\[2pt]
\overrightarrow{\widetilde{w}}(t,x)= \overrightarrow{0},\quad&  (t,x)\in \partial\Omega_T ,\\[2pt]
\overrightarrow{\widetilde{w}}(0,x)= \partial_t \overrightarrow{\widetilde{w}}(0,x)=\overrightarrow{0}, \quad& x\in \Omega,
\end{array}\right.
\end{align}
Since $\mathbb{Q}^1=\mathbb{Q}^2$ in $\mathcal{U}$, $\Xi  q_{ij}^{12}=q_{ij}^{12}$ in $\overline{\O}$ for every $i,j\in \{1,2\}$.
Also, we have $- (\Delta \Xi ) \overrightarrow{\widetilde{w}}-2 \nabla  \Xi \cdot\nabla \overrightarrow{\widetilde{w}}\equiv \overrightarrow{0}$ in $(0,T)\times (\mathcal{U}_3\cup (\Omega\setminus\overline{\mathcal{U}_2}))$, due to the construction of $\Xi $. As a result, $\overrightarrow{F}:= -\Xi \overrightarrow{b}^{12}- (\Delta \Xi ) \overrightarrow{\widetilde{w}}-2 \nabla  \Xi \cdot\nabla \overrightarrow{\widetilde{w}}$ is supported in $(0,T)\times (\O\setminus \mathcal{U})$. 
From $\overrightarrow{\widetilde{w}}\equiv\overrightarrow{0}$ in $\mathcal{U}_3$, we have 
\(
\partial_\nu \overrightarrow{\widetilde{w}}|_{\partial\O_T}=\overrightarrow{0}.
\) Now using Corollary \ref{Corollary}, we have $\overrightarrow{\widetilde{w}}=\overrightarrow{0}\text{ in }((0,T^{'})\times (\mathcal{U}_2\setminus \mathcal{U}_3)$. Next, we multiply \eqref{PDE: tilde w} with 
$\overrightarrow{w_0}$, where $\overrightarrow{w_0}$ solves 
\begin{align}
\begin{cases}
    \Box \overrightarrow{w_0}+\mathbb{P}^\top \overrightarrow{w_0}=\overrightarrow{0},& (t,x) \in \O_{T^{'}}\\ 
    \overrightarrow{w_0}(T^{'},x)=\partial_t \overrightarrow{w_0}(T^{'},x)=\overrightarrow{0}, & \quad x \in \O,
    \end{cases}
\end{align}
and integrate over $\O_{T^{'}}$ to get
\begin{align}
     \int_{\O_{T^{'}}} 
    \overrightarrow{b}^{12}\cdot \overrightarrow{w_0}~dxdt=\overrightarrow{0}.
\end{align}
Next, we choose (see Proposition \ref{Asymptotic solutions for IBVP} and \ref{Asymptotic solutions for backward IBVP})
\begin{equation}\label{eqn: IVP geom opt sol 1 in nonlinear id partial}
 \overrightarrow{v}^{(i)}(t,x)= \exp{(\mathrm{i}\lambda \phi)}\bigg(z_{0}(t,x)\overrightarrow{K}^{+,i}+\frac{Z_{1,\overrightarrow{K}^{+,i}}(t,x)}{\lambda}+\frac{Z_{2,\overrightarrow{K}^{+,i}}(t,x)}{\lambda^2}+\cdots+\frac{Z_{N,\overrightarrow{K}^{+,i}}(t,x)}{\lambda^N}\bigg)+\overrightarrow{\mathcal{R}}_{\overrightarrow{K}^{+,i}},
    \end{equation}
 and 
  \begin{equation}\label{eqn: ansatz2 nonlinear id partial}
 \overrightarrow{w_0}(t,x)= \exp{(-3\mathrm{i}\lambda \phi)}\bigg(w_{0}(t,x)\overrightarrow{K}^{-}+\frac{W_{1,\overrightarrow{K}^{-}}(t,x)}{\lambda}+\frac{W_{2,\overrightarrow{K}^{-}}(t,x)}{\lambda^2}+\cdots+\frac{W_{N,\overrightarrow{K}^{-}}(t,x)}{\lambda^N}\bigg)+\overrightarrow{\mathcal{R}}_{\overrightarrow{K}^{-}},
    \end{equation}
 where $\overrightarrow{v}^{(i)}$ and $\overrightarrow{w_0}$ are the asymptotic solutions to \eqref{PDE: v power i} and \eqref{PDE: tilde w}, respectively. Here $\varphi \in \mathtt{C}_c^\infty(\mathbb{R}^\ell)$ be arbitrary, $z_0=w_0=\varphi(x+t\theta)$ solving the transport equation \eqref{IVP: transport equation z_0} and \eqref{IVP: transport w_0}. The remainder terms  $\overrightarrow{\mathcal{R}}_{\overrightarrow{K}^{+,i}}(t,x,\lambda)$ for $1\leq i\leq 3$ and $\overrightarrow{\mathcal{R}}_{\overrightarrow{K}^{-}}(t,x,\lambda)$, satisfies the following estimates: \begin{align}\label{Bound: correction terms2}\lVert \overrightarrow{\mathcal{R}}_{\overrightarrow{K}^{+,i}}\rVert_{\textbf{L}^{\infty}(\O_T)}\leq C\lambda^{-1} \quad \mbox{ and}\quad \lVert \overrightarrow{\mathcal{R}}_{\overrightarrow{K}^{-}}\rVert_{\textbf{L}^{\infty}(\O_{T^{'}})}\leq C\lambda^{-1},\end{align} 
   where $C>0$ is independent of $\lambda$. Also, we have 
    \begin{align}  \overrightarrow{\mathcal{R}}_{\overrightarrow{K}^{+,i}}(0,x,\lambda)=\partial_t \overrightarrow{\mathcal{R}}_{\overrightarrow{K}^{+,i}}(0,x,\lambda)=\overrightarrow{0}  \mbox{ and }\; \overrightarrow{\mathcal{R}}_{\overrightarrow{K}^{-}}(T^{'},x,\lambda)=\partial_{t}\overrightarrow{\mathcal{R}}_{\overrightarrow{K}^{-}}(T^{'},x,\lambda)=\overrightarrow{0},\ \end{align}
    for $x\in \O$ and $1\leq i\leq 3$. From \eqref{Bound: correction terms2}, we observe that 
    \begin{align}
        \lVert \overrightarrow{\mathcal{R}}_{\overrightarrow{K}^{+,i}}\rVert_{\textbf{L}^{\infty}(\O_{T^{'}})}\leq C\lambda^{-1},
    \end{align}
for $1\leq i\leq 3$.
Proceeding as in \S \ref{Recovery of the coefficients of nonlinearity}, we obtain the unique identification of $\mathbb{Q}$, that is, $\mathbb{Q}^1 = \mathbb{Q}^2$ in $\O$.

    \appendix
\makeatletter
\def\thetheorem{\Alph{section}.\arabic{theorem}}
\setcounter{theorem}{0}
\makeatother
\section{Well-posedness of the forward Klein-Gordon system}\label{Well-posedness of the forward model}
In the following, we establish the well-posedness of the forward problem \eqref{Govern Eq} under smallness assumptions on the Dirichlet data. A similar well-posedness result for the scalar equation appears in \cite{bhardwaj2026reconstructionpotentialdampingcoefficients, Lassas2022UniquenessRA, lin2024determining}. Although the proof relies on quite standard arguments, we collect the proof here for the sake of completeness of the article. We move ahead by defining the function spaces and notation used in the analysis below. By referring \cite{evans2022partial}, for any $m \in \mathbb{N}$ and multi-index $\beta$ satisfying $ |\beta|\le m-1$, we define the function space
  \begin{align*}
     H^{m}_{0}(\Omega)
  := \left\{ u\in H^{m}(\Omega):\;
\partial^{\beta}u=0 \ \text{on}\ \partial\Omega,
            \ \text{ for } |\beta|\le m-1 \right\}, 
\end{align*}
as the completion  of $(C_c^{\infty}(\Omega), \lVert \cdot \rVert_{H^m(\Omega)})$ with $H^m$ representing the standard Sobolev space. In the case of vector-valued functions, we define the function space
  \begin{align*}
     \textbf{H}^{m}_{0}(\Omega)
  := \left\{ \overrightarrow{u}\in \textbf{H}^{m}(\Omega):\;
\partial^{\beta}\overrightarrow{u}=\overrightarrow{0} \ \text{on}\ \partial\Omega,
            \ \text{ for } |\beta|\le m-1 \right\} 
\end{align*}
as the completion of $(\textbf{C}_c^{\infty}(\Omega), \lVert \cdot \rVert_{\textbf{H}^m(\Omega)})$, where $\lVert \cdot \rVert_{\textbf{H}^m(\O)}$ is defined as  \[\lVert\overrightarrow{u}\rVert_{\textbf{H}^m(\O)}^2:=\lVert u_1\rVert_{H^m(\O)}^2+\lVert u_2\rVert_{H^m(\O)}^2.\]
Consider the linear inhomogeneous initial-boundary value problem (IBVP)
\begin{equation}\label{PDE: inhomogenuous pde for scalar case}
\begin{cases}
\displaystyle \Box u(t,x) = G(t,x),& (t,x)\in \Omega_T,  \\[2pt]
u(t,x)= f(t,x),&  (t,x)\in \partial\Omega_T ,\\[2pt]
u(0,x)= 0, \quad\partial_t u(0,x)= 0, &\quad x\in \Omega.
\end{cases}
\end{equation}
Let  $D_m$ denote the energy space, which is given by 
\begin{align*}
  D_m
    &:= \bigcap_{s=0}^m C^s\bigl([0,T];H^{m-s}(\Omega)\bigr),
\end{align*}
equipped with the norm 
\begin{align}\label{norm-Em}
\|u\|_{D_m}^2 
    := \sup_{t\in[0,T]} 
       \sum_{s=0}^{m} 
       \bigl\| \partial_t^{\,s} u(t,\cdot) \bigr\|_{H^{m-s}(\Omega)}^2.
\end{align}
By referring to  \cite[Corollary~2]{Lassas2022UniquenessRA} and \cite[Lemma~ 3.2]{lin2024determining}, we have 
the following well-posedness result stated below.
\begin{lemma}\cite[Corollary ~2]{Lassas2022UniquenessRA}
Let  
 $T>0$, $m>\ell+1$. Then for any $f\in H^{m+1}(\partial\O_T)$ with $\partial_t^k f(0,\cdot)=0$ on $\partial\O$ for $k=0,1,\cdots,m$ and $G\in D_{m+1}$ with $\partial_t^k G(0,\cdot)=0$ in $\O$ for $k=0,1,\cdots,m$, 
the IBVP \eqref{PDE: inhomogenuous pde for scalar case}
admits a unique solution $u\in D_{m+1}$ satisfying $\partial_\nu u \in H^m(\partial \O_T)$ and \begin{align}\label{energy estimate-v}
\begin{aligned}
   \|u\|_{D_{m+1}} + \|\partial_\nu u\|_{H^m(\partial \O_T)} \le C\,\exp{(CT)}\left(\sum_{k=0}^{m}\|\partial_t^{\,k} G\|_{L^1(0,T;H^{m-k}(\Omega))}
 + \|f\|_{H^{m+1}(\partial \O_T)}\right),
\end{aligned}
\end{align}
where $C>0$,  depends only on   $\Omega$.
\end{lemma}
Clearly, for a decoupled linear inhomogeneous wave system
\begin{equation}\label{PDE: inhomogenuous pde for system case}
\begin{cases}
\displaystyle \Box\overrightarrow{u}(t,x) = \overrightarrow{F}(t,x),& (t,x)\in \Omega_T,  \\[2pt]
\overrightarrow{u}(t,x)= \overrightarrow{f}(t,x),&  (t,x)\in \partial\Omega_T ,\\[2pt]
\overrightarrow{u}(0,x)= \overrightarrow{0}, \quad\partial_t\overrightarrow{u}(0,x)= \overrightarrow{0}, &\quad x\in \Omega,
\end{cases}
\end{equation}
one can extend the above result and obtain
\begin{lemma}\label{Lemma 1 system}
 Let  
 $T>0$, $m>\ell+1$. Then for any $\overrightarrow{f}\in \textbf{H}^{m+1}(\partial\O_T)$ with $\partial_t^k \overrightarrow{f}(0,\cdot)=\overrightarrow{0}$ on $\partial\O$ for $k=0,1,\cdots,m$ and $\overrightarrow{F}\in \mathscr{D}_{m+1}$ with $\partial_t^k \overrightarrow{F}(0,\cdot)=\overrightarrow{0}$ in $\O$ for $k=0,1,\cdots,m$, 
the IBVP \eqref{PDE: inhomogenuous pde for system case}
admits a unique solution  $\overrightarrow{u}\in \mathscr{D}_{m+1}$ satisfying $\partial_\nu \overrightarrow{u} \in \textbf{H}^m(\partial \O_T)$ and \begin{align}\label{energy estimate-overright u}
\begin{aligned}
   \|\overrightarrow{u}\|_{\mathscr{D}_{m+1}} + \|\partial_\nu \overrightarrow{u}\|_{\textbf{H}^m(\partial \O_T)} \le C\,\exp{(CT)}\left(\sum_{k=0}^{m}\|\partial_t^{\,k} \overrightarrow{F}\|_{L^1(0,T;\textbf{H}^{m-k}(\Omega))}
+ \|\overrightarrow{f}\|_{\textbf{H}^{m+1}(\partial \O_T)}\right),
\end{aligned}
\end{align}
where 
$C$ is a positive constant,  depends only on  $\Omega$ and energy space $\mathscr{D}_m$ is given by 
 \begin{align}\label{Energy space}
      \mathscr{D}_m
    := \bigcap_{s=0}^m C^s\bigl([0,T];\textbf{H}^{m-s}(\Omega)\bigr),
 \end{align}
 endowed with the norm
 \begin{align}\label{norm-Em_1}
\|\overrightarrow{u}\|_{\mathscr{D}_m}^2
    := \sup_{t\in[0,T]} 
       \sum_{s=0}^{m} \sum_{i=1}^2
       \bigl\| \partial_t^{\,s} u_i(t,\cdot) \bigr\|_{H^{m-s}(\Omega)}^2.
\end{align}
\end{lemma}
Next, we introduce a zeroth-order perturbation to \eqref{PDE: inhomogenuous pde for system case} and seek the well-posedness of the IBVP 
\begin{align}\label{PDE: inhomogenuous pde for system case with zero order perturbation}
\left\{\begin{array}{r l c}
\displaystyle \mathcal{L}_{\mathbb{P}}\overrightarrow{u}(t,x)& = \overrightarrow{F}(t,x),& (t,x)\in \Omega_T,  \\[2pt]
\overrightarrow{u}(t,x)&= \overrightarrow{f}(t,x),&  (t,x)\in \partial\Omega_T ,\\[2pt]
\overrightarrow{u}(0,x)&= \partial_t\overrightarrow{u}(0,x)= \overrightarrow{0}, &\quad x\in \Omega,
\end{array}\right.
\end{align}
which is established in the proposition given below.
\begin{proposition}\label{proposition 1 system}
In addition to notations and assumption of Lemma \ref{Lemma 1 system}, we further assume that $ \mathbb{P}=\{p_{ij}\}_{1\leq i,j\leq 2}\in \textbf{C}^{\infty}_c(\O)$. Then the IBVP \eqref{PDE: inhomogenuous pde for system case with zero order perturbation}
admits a unique solution  $\overrightarrow{u}\in \mathscr{D}_{m+1}$ satisfying $\partial_\nu \overrightarrow{u} \in \textbf{H}^m(\partial \O_T)$ and \begin{align}\label{energy estimate-overright u 1}
\begin{aligned}
   \|\overrightarrow{u}\|_{\mathscr{D}_{m+1}} + \|\partial_\nu \overrightarrow{u}\|_{\textbf{H}^m(\partial \O_T)} \le C\,\exp{(CT)}\left(\sum_{k=0}^{m}\|\partial_t^{\,k} \overrightarrow{F}\|_{L^1(0,T;\textbf{H}^{m-k}(\Omega))} + \|\overrightarrow{f}\|_{\textbf{H}^{m+1}(\partial \O_T)}\right),
\end{aligned}
\end{align}
where $C>0$,  depends only on  $\Omega$ and $\mathbb{P}$.
\end{proposition}
\begin{proof}
For a given $\overrightarrow{v}\in \mathscr{D}_{m+1}$ with $\partial_t^k \overrightarrow{v}(0,\cdot) =0$ in $\O$ for $k=0,1,\cdots,m$
, we consider the IBVP 
\begin{equation}\label{PDE: IBVP with rhs}
\begin{cases}
\displaystyle \Box\overrightarrow{u}(t,x) = \overrightarrow{F}(t,x)-\mathbb{P}(x)\overrightarrow{v}(t,x),& (t,x)\in \Omega_T,  \\[2pt]
\overrightarrow{u}(t,x)= \overrightarrow{f}(t,x),&  (t,x)\in \partial\Omega_T ,\\[2pt]
\overrightarrow{u}(0,x)= \overrightarrow{0}, \quad\partial_t\overrightarrow{u}(0,x)= \overrightarrow{0}, &\quad x\in \Omega.
\end{cases}
\end{equation}
The idea is to use the contraction mapping theorem to establish local existence of the solution. We can obtain global existence by rescaling with respect to the time variable.
For $m>\ell+1$, $\mathscr{D}_{m+1}$ is an Banach algebra (see \cite[Definition 3.5]{choquet2009general}). Given $\mathbb{P}\in \textbf{C}^{\infty}_c(\O)$ and $\overrightarrow{F}\in \mathscr{D}_{m+1}$ with $\partial_t^k \overrightarrow{F}(0,\cdot)=\overrightarrow{0}$ in $\O$ for $k=0,1,\cdots,m$, we have $\overrightarrow{F}-\mathbb{P}\overrightarrow{v}\in \mathscr{D}_{m+1}$. Although, to apply Lemma \ref{Lemma 1 system}, we further need to establish that $\partial_t^k( \overrightarrow{F}(0,\cdot)-\mathbb{P}\overrightarrow{v}(0,\cdot))=\overrightarrow{0}$ in $\O$ for $k=0,1,\cdots,m$, which holds trivially due to choice of  $\overrightarrow{v}$.  Thus the inhomogeneous IBVP \eqref{PDE: IBVP with rhs} admits a unique solution $\overrightarrow{u}\in \mathscr{D}_{m+1}$ with $\partial_\nu \overrightarrow{u} \in \textbf{H}^m(\partial \O_T)$ and
\begin{align}\label{Energy estiomate}
\begin{aligned}
\|\overrightarrow{u}\|_{\mathscr{D}_{m+1}} &+ \|\partial_\nu \overrightarrow{u}\|_{\textbf{H}^m(\partial \O_T)}\\& \le C\,\exp{(CT)}\Big(\sum_{k=0}^{m}\|\partial_t^{\,k}( \overrightarrow{F}-\mathbb{P} \overrightarrow{v})\|_{L^1(0,T;\textbf{H}^{m-k}(\Omega))}+ \|\overrightarrow{f}\|_{\textbf{H}^{m+1}(\partial \O_T)}\Big).
\end{aligned}
\end{align}
Set 
\[\mathscr{N}= \{\overrightarrow{v}\in  \mathscr{D}_{m+1} : \partial_t^k \overrightarrow{v}(0,\cdot) =0,~k=0,1,\cdots,m\}.\]
Define a map
\(
    \mathscr{T}: \mathscr{N}\rightarrow\mathscr{N}
\)
mapping
\begin{align*}
    \mathscr{T}(\overrightarrow{v})=\overrightarrow{u},
\end{align*}
where $\overrightarrow{u}$ is the solution to the IBVP \eqref{PDE: IBVP with rhs} corresponds to $\overrightarrow{v}\in \mathscr{N}$. Although $\overrightarrow{u}\in \mathscr{D}_{m+1}$, the compatibility condition of \eqref{PDE: IBVP with rhs} forces $\overrightarrow{u}$ to lie in $\mathscr{N}$, hence, $\mathscr{T}$ is well-defined.
As a next step, we establish that $\mathscr{T}$ is a contraction map. To this end, let $\overrightarrow{u}^i$ be the solution corresponds to the IBVP  
\begin{equation}\label{PDE: IBVP with rhs 1}
\begin{cases}
\displaystyle \Box\overrightarrow{u}^i(t,x) = \overrightarrow{F}(t,x)-\mathbb{P}(x)\overrightarrow{v}^i(t,x),& (t,x)\in \Omega_T,  \\[2pt]
\overrightarrow{u}^i(t,x)= \overrightarrow{f}(t,x),&  (t,x)\in \partial\Omega_T ,\\[2pt]
\overrightarrow{u}^i(0,x)= \overrightarrow{0}, \quad\partial_t\overrightarrow{u}^i(0,x)= \overrightarrow{0}, &\quad x\in \Omega,
\end{cases}
\end{equation}
for $i=1,2,$ and $\overrightarrow{v}^i \in \mathscr{N}$. The subtraction then leads us to
\begin{equation}
\begin{cases}
\displaystyle \Box(\overrightarrow{u}^1-\overrightarrow{u}^2)(t,x) =-\mathbb{P}(x)(\overrightarrow{v}^1-\overrightarrow{v}^2)(t,x),& (t,x)\in \Omega_T,  \\[2pt]
(\overrightarrow{u}^1-\overrightarrow{u}^2)(t,x)= \overrightarrow{0},&  (t,x)\in \partial\Omega_T ,\\[2pt]
(\overrightarrow{u}^1-\overrightarrow{u}^2)(0,x)= \overrightarrow{0}, \quad\partial_t(\overrightarrow{u}^1-\overrightarrow{u}^2)(0,x)= \overrightarrow{0}, &\quad x\in \Omega,
\end{cases}
\end{equation}
and the energy estimate
\begin{align}\label{Energy estimate2}
\begin{aligned}
\|\overrightarrow{u}^1-\overrightarrow{u}^2\|_{\mathscr{D}_{m+1}} + \|\partial_\nu (\overrightarrow{u}^1-\overrightarrow{u}^2)\|_{\textbf{H}^m(\partial \O_T)}
    &\le C\,\exp{(CT)}\Big(\sum_{k=0}^{m}\|\mathbb{P}\partial_t^{\,k}( \overrightarrow{v}^1-\overrightarrow{v}^2)\|_{L^1(0,T;\textbf{H}^{m-k}(\Omega))}
\Big)\\&\leq C T\,\exp{(CT)}\Big(\sum_{k=0}^{m}\|\mathbb{P}\partial_t^{\,k}( \overrightarrow{v}^1-\overrightarrow{v}^2)\|_{C(0,T;\textbf{H}^{m-k}(\Omega))}
\Big)\\&\leq CC_1 T\,\exp{(CT)} ||\mathbb{P}||_{\textbf{H}^m(\O)}||(\overrightarrow{v}^1-\overrightarrow{v}^2)||_{\mathscr{D}_{m+1}},
\end{aligned}
\end{align}
by utilizing the estimate $\|\overrightarrow{w}\|_{L^1(0,T;\textbf{H}^{m-k}(\Omega))}\leq T \|\overrightarrow{w}\|_{C(0,T;\textbf{H}^{m-k}(\Omega))}$.
For a sufficient small $T$, say $T=T^*$, we have $C C_1T^* ||\mathbb{P}||_{\textbf{H}^m(\O)}\exp{(CT^*)}<1$, resulting in 
\begin{align}
  \|\overrightarrow{u}^1-\overrightarrow{u}^2\|_{\mathscr{D}_{m+1}([0, T^*])} <   ||(\overrightarrow{v}^1-\overrightarrow{v}^2)||_{\mathscr{D}_{m+1}( [0, T^*])},
\end{align}
where $[0, T^*]$ with $\mathscr{D}_m$ indicated that the initial and final time are $0$ and $T^*$ respectively in the aforementioned definition of  $\mathscr{D}_m$.
Thus $\mathscr{T}$ is a contraction map, and by using the Banach fixed point theorem, there exists a unique fixed point $\overrightarrow{v}^{(0)}\in \mathscr{N}$, that is $\mathscr{T}(\overrightarrow{v}^{(0)})=\overrightarrow{v}^{(0)}=\overrightarrow{u}$. As a result, \eqref{PDE: inhomogenuous pde for system case with zero order perturbation} admits a unique solution whenever $0\leq t\leq T^*$. Moreover, using the fixed point $\overrightarrow{v}$, and the estimate \eqref{Energy estiomate}, we have
\begin{align}\label{Energy estimate 1}
\begin{aligned}
&\|\overrightarrow{u}\|_{\mathscr{D}_{m+1}([0, T^*])} + \|\partial_\nu \overrightarrow{u}\|_{\textbf{H}^m(\partial \O_{T^*})}\\
   & \le C\,\exp{(CT^*)}\Big(\sum_{k=0}^{m}\|\partial_t^{\,k}( \overrightarrow{F}-\mathbb{P} \overrightarrow{u})\|_{L^1(0,T^*;\textbf{H}^{m-k}(\Omega))}
 + \|\overrightarrow{f}\|_{\textbf{H}^{m+1}(\partial \O_{T^{*}})}\Big)\\&
\le C\,\exp{(CT^*)}\Big(\sum_{k=0}^{m}\|\partial_t^{\,k} \overrightarrow{F}\|_{L^1(0,T^*;\textbf{H}^{m-k}(\Omega))}
+  C_1T^* ||\mathbb{P}||_{\textbf{H}^m(\O)}||\overrightarrow{u}||_{\mathscr{D}_{m+1}( [0, T^*])} + \|\overrightarrow{f}\|_{\textbf{H}^{m+1}(\partial \O_{T^{*}})}\Big).
\end{aligned}
\end{align}
From the choice of $T^*$, the aforementioned inequality boils down to
\begin{align}
\begin{aligned}
\|\overrightarrow{u}\|_{\mathscr{D}_{m+1}( [0, T^*])} + \|\partial_\nu \overrightarrow{u}\|_{\textbf{H}^m(\partial \O_{T^*})}
\le C\,\exp{(CT^*)}\Big(\sum_{k=0}^{m}\|\partial_t^{\,k} \overrightarrow{F}\|_{L^1(0,T^*;\textbf{H}^{m-k}(\Omega))} + \|\overrightarrow{f}\|_{\textbf{H}^{m+1}(\partial \O_{T^{*}})}\Big).
\end{aligned}
\end{align}
Since the IBVP \eqref{PDE: inhomogenuous pde for system case with zero order perturbation} is linear, rescaling in the time variable (see \cite[Lemma~~3.2]{bhardwaj2026reconstructionpotentialdampingcoefficients}) completes the proof.
\end{proof}
Next, we introduce the small  Dirichlet data $\overrightarrow{f}\in \mathtt{F}^\sigma_{m+1}$, where
\begin{align}\label{E-delta}
\begin{aligned}
\mathtt{F}^\sigma_{m+1}
:= \bigg\{ \overrightarrow{f}\in \textbf{H}^{m+1}(\partial \O_T)\;:\; 
\partial_t^k \overrightarrow{f}(0,\cdot)=\overrightarrow{0} \text{ on } \partial\O  \text{ and }\ \|\overrightarrow{f}\|_{\textbf{H}^{m+1}(\partial\Omega_T)}< \sigma,\\
\mbox{for } k=0,1,\dots,m \bigg\},
\end{aligned}
\end{align}
with $\sigma>0$ chosen sufficiently small. Further, we denote
\begin{align*}
    \mathbb{B}_R(0) &:= \left\{ \overrightarrow{u} \in \mathscr{D}_{m+1} \;:\; \|\overrightarrow{u}\|_{\mathscr{D}_{m+1}} <R\right\} \text{ for some } R>0.
\end{align*}
Using the preceding analysis, we establish a unique solution $\overrightarrow{u}\in     \mathbb{B}_R(0)$ for the IBVP  \eqref{Govern Eq},
whenever $\overrightarrow{f}\in \mathtt{F}^\sigma_{m+1}$, which is stated formally below.
\begin{theorem}\label{Thm: forward main wellposedness}
    Let $T>0$, $m>\ell+1$, $(\mathbb{P}, \mathbb{Q})\in \textbf{C}^{\infty}_c(\O)\times \textbf{C}^{\infty}_c(\O)$, and $\overrightarrow{f}\in \mathtt{F}^\sigma_{m+1}$.
Then the IBVP 
\eqref{Govern Eq}
 admits a unique solution $\overrightarrow{u}\in  \mathbb{B}_R(0)$ satisfying $\partial_\nu \overrightarrow{u} \in \textbf{H}^m(\partial \O_T)$ and \begin{align}\label{energy estimate main}
\begin{aligned}
   \|\overrightarrow{u}\|_{\mathscr{D}_{m+1}} + \|\partial_\nu \overrightarrow{u}\|_{\textbf{H}^m(\partial \O_T)}\le C\,\exp{(CT)} \|\overrightarrow{f}\|_{\textbf{H}^{m+1}(\partial \O_T)},
\end{aligned}
\end{align}
where $C>0$,  depends only on  $\Omega, P$ and $Q$.
\end{theorem}
\begin{proof}
We begin by introducing a closed subset of $\mathscr{D}_{m+1}$ which is given by
\begin{align}
    \mathscr{K}:=\Big\{\overrightarrow{v} \in \mathbb{B}_R(0):  \partial_t^k \overrightarrow{v}(0,\cdot) =0,~k=0,1,\cdots,m\Big\},
\end{align}
where $R>0$ is specified later.
Given $\overrightarrow{v}\in   \mathscr{K}$, we consider the IBVP
\begin{equation}\label{PDE: inhomogenuous nonlinear pde for system case with zero order perturbation}
\begin{cases}
\displaystyle \Box\overrightarrow{u}(t,x)+\mathbb{P}(x)\overrightarrow{u}(t,x) = -(\mathbb{Q}\overrightarrow{v})\circ(\overrightarrow{v}\circ \overrightarrow{v})(t,x),& (t,x)\in \Omega_T,  \\[2pt]
\overrightarrow{u}(t,x)= \overrightarrow{f}(t,x),&  (t,x)\in \partial\Omega_T ,\\[2pt]
\overrightarrow{u}(0,x)= \overrightarrow{0}, \quad\partial_t\overrightarrow{u}(0,x)= \overrightarrow{0}, &\quad x\in \Omega.
\end{cases}
\end{equation}
Since $\overrightarrow{v}\in   \mathscr{K}$ and $\mathbb{P}, \mathbb{Q}\in \textbf{C}_c^\infty(\O)$, it follows from Proposition \ref{proposition 1 system} that the IBVP \eqref{PDE: inhomogenuous nonlinear pde for system case with zero order perturbation} admits a unique solution $\overrightarrow{u}\in \mathscr{D}_{m+1}$ with $\partial_\nu \overrightarrow{u} \in \textbf{H}^m(\partial \O_T)$ and
\begin{align}\label{Energy estiomate 7}
\begin{aligned}
\|\overrightarrow{u}\|_{\mathscr{D}_{m+1}}+ \|\partial_\nu \overrightarrow{u}\|_{\textbf{H}^m(\partial \O_T)} \le C\,\exp{(CT)}\Big(\sum_{k=0}^{m}\|\partial_t^{\,k}(\mathbb{Q}\overrightarrow{v})\circ(\overrightarrow{v}\circ \overrightarrow{v})\|_{L^1(0,T;\textbf{H}^{m-k}(\Omega))}
+ \|\overrightarrow{f}\|_{\textbf{H}^{m+1}(\partial \O_T)}\Big),
\end{aligned}
\end{align}
where $C>0$, depends only on $\O,\mathbb{P}$ and $\mathbb{Q}$. As discussed earlier, we now define a map \(
    \mathscr{T}: \mathscr{K}\rightarrow\mathscr{K}
\), which is a contraction. For a given $\overrightarrow{v}\in \mathscr{K}$, let $ \mathscr{T} \overrightarrow{v}=\overrightarrow{u}$, where $\overrightarrow{u} $ is a solution to the IBVP \eqref{PDE: inhomogenuous nonlinear pde for system case with zero order perturbation}. However, the well-posedness result implies that $\overrightarrow{u} \in \mathscr{D}_{m+1}$. So, to ensure that the map $\mathscr{T}$ is a contraction from $\mathscr{K}$ to itself, we need to impose some conditions on $R$ and $\sigma.$ From the above energy estimate, we have 
\begin{align}\label{Energy estiomate 8}
\begin{aligned}
\|\mathscr{T}\overrightarrow{v}\|_{\mathscr{D}_{m+1}}=\|\overrightarrow{u}\|_{\mathscr{D}_{m+1}}  \le C\,\exp{(CT)}\Big(C_1 R^3||\mathbb{Q}||_{\textbf{H}^m(\O)}
 + \sigma \Big).
\end{aligned}
\end{align}
To ensure $\overrightarrow{u}\in \mathscr{K}$, we need to have $\|\overrightarrow{u}\|_{\mathscr{D}_{m+1}} \leq R$. To make this happen, we choose 
\begin{align}\label{condition: R and sigma}
    R<\frac{1}{\sqrt{3CC_1 ||\mathbb{Q}||_{\textbf{H}^m(\O)}\exp{(CT)} }} \text{ and } \sigma<\frac{R}{3C \exp{(CT)}}.
\end{align}
Using the above choices of $R$ and $\sigma$ in \eqref{Energy estiomate 8} leads to
\begin{align}
    \|\mathscr{T}\overrightarrow{v}\|_{\mathscr{D}_{m+1}}\leq \frac{2}{3}R<R.
\end{align}
From now on, we will work only with the aforementioned choice of $R$ and $\sigma$.
Next, we set $\overrightarrow{u}^i: =\mathscr{T}\overrightarrow{v}^i$ for $\overrightarrow{v}^i\in \mathscr{K}$ and let $\overrightarrow{u}^i$ be the solution corresponds to the IBVP  
\begin{equation}\label{PDE: IBVP with rhs 2}
\begin{cases}
\displaystyle \Box\overrightarrow{u}^i(t,x)+\mathbb{P}(x)\overrightarrow{u}^i(t,x) = -(\mathbb{Q}\overrightarrow{v}^i)\circ(\overrightarrow{v}^i\circ \overrightarrow{v}^i)(t,x),& (t,x)\in \Omega_T,  \\[2pt]
\overrightarrow{u}^i(t,x)= \overrightarrow{f}(t,x),&  (t,x)\in \partial\Omega_T ,\\[2pt]
\overrightarrow{u}^i(0,x)= \overrightarrow{0}, \quad\partial_t\overrightarrow{u}^i(0,x)= \overrightarrow{0}, &\quad x\in \Omega,
\end{cases}
\end{equation}
for $i=1,2,$. The subtraction then leads us to
\begin{equation}
\begin{cases}
\displaystyle \Box(\overrightarrow{u}^1-\overrightarrow{u}^2) +\mathbb{P}(\overrightarrow{u}^1-\overrightarrow{u}^2)=-(\mathbb{Q}\overrightarrow{v}^1)\circ(\overrightarrow{v}^1\circ \overrightarrow{v}^1)+ (\mathbb{Q}\overrightarrow{v}^2)\circ(\overrightarrow{v}^2\circ \overrightarrow{v}^2),& (t,x)\in \Omega_T,  \\[2pt]
(\overrightarrow{u}^1-\overrightarrow{u}^2)(t,x)= \overrightarrow{0},&  (t,x)\in \partial\Omega_T ,\\[2pt]
(\overrightarrow{u}^1-\overrightarrow{u}^2)(0,x)= \overrightarrow{0}, \quad\partial_t(\overrightarrow{u}^1-\overrightarrow{u}^2)(0,x)= \overrightarrow{0}, &\quad x\in \Omega,
\end{cases}
\end{equation}
and the energy estimate 
\begin{align}\label{Energy estimate3}
\begin{aligned}
&\|\overrightarrow{u}^1-\overrightarrow{u}^2\|_{\mathscr{D}_{m+1}} + \|\partial_\nu (\overrightarrow{u}^1-\overrightarrow{u}^2)\|_{\textbf{H}^m(\partial \O_T)}
    \\&\le C\,\exp{(CT)}\Big(\sum_{k=0}^{m}\|\partial_t^{\,k}\left[ -(\mathbb{Q}\overrightarrow{v}^1)\circ(\overrightarrow{v}^1\circ \overrightarrow{v}^1)+ (\mathbb{Q}\overrightarrow{v}^2)\circ(\overrightarrow{v}^2\circ \overrightarrow{v}^2)\right]\|_{L^1(0,T;\textbf{H}^{m-k}(\Omega))}
\Big).
\end{aligned}
\end{align}
If we denote $\overrightarrow{v}^i:=(v^i_1,v^i_2)^\top$ then
\begin{align}
    &-(\mathbb{Q}\overrightarrow{v}^1)\circ(\overrightarrow{v}^1\circ \overrightarrow{v}^1)+ (\mathbb{Q}\overrightarrow{v}^2)\circ(\overrightarrow{v}^2\circ \overrightarrow{v}^2)\\[4pt]&=-\begin{pmatrix}
        q_{11}((v^1_1)^3-(v^2_1)^3)+ q_{12}((v_1^1)^2(v_2^1)-(v_1^2)^2(v_2^2))\\[3pt]q_{21}((v^1_2)^2(v_1^1)-(v_2^2)^2(v_1^2))+q_{22}((v_2^1)^3-(v_2^2)^3)
    \end{pmatrix}\\[4pt]&=-\begin{pmatrix}
        q_{11}(v_1^1-v_1^2)((v_1^1)^2+(v^2_1)^2+v_1^1v_1^2)+ q_{12}(v_1^1-v_1^2)(v_1^1+v_1^2)v_2^1+q_{12}(v_2^1-v_2^2)(v_1^2)^2\\[4pt]  q_{21}(v_1^1-v_1^2)(v^1_2)^2+q_{21}(v^2_1)(v_2^1-v_2^2)(v_2^1+v_2^2)+q_{22}(v_2^1-v_2^2)((v_2^1)^2+(v_2^2)^2+v_2^1v_2^2)
    \end{pmatrix}.
\end{align}
The use of foregoing expression along with the fact that $\overrightarrow{v}^i \in \mathscr{K}$ in \eqref{Energy estimate3} leads to
\begin{align}
\begin{aligned}
&\lVert\mathscr{T}\overrightarrow{v}^1-\mathscr{T}\overrightarrow{v}^2\rVert_{\mathscr{D}_{m+1}}= \|\overrightarrow{u}^1-\overrightarrow{u}^2\|_{\mathscr{D}_{m+1}}
    \le CC_1 R^2||\mathbb{Q}||_{\textbf{H}^m(\O)}\exp{(CT)}\lVert \overrightarrow{v}^1-\overrightarrow{v}^2\rVert_{\mathscr{D}_{m+1}}
 .
\end{aligned}
\end{align}
If we choose 
\begin{align}
    R< \dfrac{1}{2\sqrt{CC_1 ||\mathbb{Q}||_{\textbf{H}^m(\O)}\exp{(CT)}}},
\end{align}
then
\begin{align}
\begin{aligned}
&\lVert\mathscr{T}\overrightarrow{v}^1-\mathscr{T}\overrightarrow{v}^2\rVert_{\mathscr{D}_{m+1}}
< \lVert \overrightarrow{v}^1-\overrightarrow{v}^2\rVert_{\mathscr{D}_{m+1}}
 .
\end{aligned}
\end{align}
Thus, the Banach fixed point theorem ensures the existence of a unique $\overrightarrow{v}\in \mathscr{K}$ such that $\mathscr{T}(\overrightarrow{v})=\overrightarrow{v}$. Thus there exist a unique  $\overrightarrow{u}\in \mathscr{D}_{m+1}$ solving \eqref{Govern Eq}, such that
\begin{align}
\begin{aligned}
\|\overrightarrow{u}\|_{\mathscr{D}_{m+1}}+ \|\partial_\nu \overrightarrow{u}\|_{\textbf{H}^m(\partial \O_T)} &\le C\,\exp{(CT)}\Big(\sum_{k=0}^{m}\|\partial_t^{\,k}(\mathbb{Q}\overrightarrow{u})\circ(\overrightarrow{u}\circ \overrightarrow{u})\|_{L^1(0,T;\textbf{H}^{m-k}(\Omega))}
+ \|\overrightarrow{f}\|_{\textbf{H}^{m+1}(\partial \O_T)}\Big)\\&\leq C\,\exp{(CT)}\|(\mathbb{Q}\overrightarrow{u})\circ(\overrightarrow{u}\circ \overrightarrow{u})\|_{\mathscr{D}_{m}}+C\,\exp{(CT)}\|\overrightarrow{f}\|_{\textbf{H}^{m+1}(\partial \O_T)}\\&\leq  2CC_1R^2||\mathbb{Q}||_{\textbf{H}^m(\O)}\exp{(CT)}\lVert \overrightarrow{u}\rVert_{\mathscr{D}_{m+1}}+C\,\exp{(CT)}\|\overrightarrow{f}\|_{\textbf{H}^{m+1}(\partial \O_T)}\\& \leq \frac{1}{2}\lVert \overrightarrow{u}\rVert_{\mathscr{D}_{m+1}}+C\,\exp{(CT)}\|\overrightarrow{f}\|_{\textbf{H}^{m+1}(\partial \O_T)}.
\end{aligned}
\end{align}
Thus, we have  
\begin{align}\label{Energy estiomate4}
\begin{aligned}
\|\overrightarrow{u}\|_{\mathscr{D}_{m+1}}+ \|\partial_\nu \overrightarrow{u}\|_{\textbf{H}^m(\partial \O_T)} \leq C\,\exp{(CT)}\|\overrightarrow{f}\|_{\textbf{H}^{m+1}(\partial \O_T)}.
\end{aligned}
\end{align}
This completes the proof.
\end{proof}
\small
\bibliographystyle{plain}
\bibliography{references}
\end{document}